\documentclass[12pt, reqno, draft]{amsart}

\usepackage{ifdraft}
\usepackage{fullpage}

\usepackage{graphicx}
\usepackage[dvipsnames]{xcolor}
\usepackage{pgf}

 \usepackage{tikz}
 \usetikzlibrary{matrix,arrows,decorations.pathmorphing}
 \usepackage{tikz-cd}
 \usetikzlibrary{arrows}

\usepackage{comment}

\usepackage{xypic}

\newcommand{\cat}[1]{\mathbf{#1}}

\usepackage{amsmath}
\usepackage{mathrsfs}

\usepackage{epsfig}
\usepackage{amsfonts,amsthm,amscd, amssymb}
\usepackage{latexsym}
\usepackage{enumerate}
\usepackage{mypack}

\newcommand{\uHom}{\underline{\operatorname{Hom}}}
\newcommand{\agroupoid}{/\!\!/}
\newcommand{\iso}{\cong}

\usepackage[colorinlistoftodos,prependcaption,textsize=tiny]{todonotes} 

\newcommand{\benw}[2][]{\ifdraft{\todo[linecolor=Green,backgroundcolor=Green!25,bordercolor=Green,#1]{#2---Ben W.}}{}}

\theoremstyle{definition}
\newtheorem{example}[thm]{Example}

\numberwithin{equation}{thm}
\title{On the mod-$\ell$ homology of the classifying space for commutativity}
\date{\today}
\address{Department of Physics and Astronomy, University of British Columbia, Vancouver BC V6T 1Z4, Canada}
\email{okay@math.ubc.ca/okay@phas.ubc.ca}
\author{C\.{I}han Okay}
\author{Ben Williams}
\address{Department of Mathematics, University of British Columbia, Vancouver BC V6T 1Z2, Canada}
\email{tbjw@math.ubc.ca}

\begin{document}

\begin{abstract}
   We study the mod-$\ell$ homotopy type of classifying spaces for commutativity, $B(\ZZ, G)$, at a prime $\ell$. We show that the
   mod-$\ell$ homology of $B(\ZZ, G)$ depends on the mod-$\ell$ homotopy type of $BG$ when $G$ is a compact connected
   Lie group, in the sense that a mod-$\ell$ homology isomorphism $BG \to BH$ for such groups induces a mod-$\ell$ homology isomorphism
   $B(\ZZ, G) \to B(\ZZ, H)$. In order to prove this result, we study a presentation of $B(\ZZ, G)$ as a homotopy
   colimit over a topological poset of closed abelian subgroups, expanding on an idea of Adem and G\'omez. We also study
   the relationship between the mod-$\ell$ type of a Lie group $G(\CC)$ and the locally finite group
   $G(\bar{\mathbb{F}}_p)$ where $G$ is a Chevalley group. We see that the na\"ive analogue for $B(\ZZ, G)$ of the
   celebrated Friedlander--Mislin result cannot hold, but we show that it does hold after taking the homotopy quotient
   of a $G$ action on $B(\ZZ, G)$.
 \end{abstract}

 \maketitle

\section{Introduction }
 Suppose $G$ is a topological group. Following \cite{AG15}, it is possible to define a space $B_{\text{com}}(G)$ or
 $B(\ZZ, G)$, a realization of a bar construction after the fashion of $BG$, but where the space of $n$-simplices
 consists only of $n$-tuples of pairwise commuting elements of $G$. This construction and its variations are first
 introduced in \cite{ACT12}, and further studied in \cite{AG15,AGLT,Gri17,AGV17} when $G$ is a compact (connected) Lie
 group; whereas the finite group case is considered in \cite{O14,O15,O16}.

The construction of $B(\ZZ,G)$ is part of a family of constructions: given a topological group $G$ and a cosimplicial
group $\tau^\dt$, finitely generated in each degree, one may form a simplicial space $B(\tau, G)_\dt$, having
$\Hom(\tau^n, G)$ in the $n$-th degree, and then realize to form $B(\tau, G)$. There exists a cosimplicial group
$F^\dt$, free of rank $n$ in degree $n$, for which $B(F,G)$ is the usual bar construction of $BG$. Replacing $F^\dt$ by
its abelianization yields $\ZZ^\dt$ and $B(\ZZ,G)$. For any integer $m$, one may reduce further to produce
$B(\ZZ/m, G)$, which plays an important auxiliary role in this paper. That is, this   analogue of $BG$ one obtains
by restricting attention to $n$-tuples of elements of $G$, each of which commutes with each of the others and each of
which has order dividing $m$.

This paper arose from our attempt to understand the mod-$\ell$ homology of $B(\ZZ, G)$ for a prime $\ell$ where $G$ is a compact Lie group. We
establish that when $\pi_0(G)$ is an $\ell$-group, the mod-$\ell$ homotopy type of $B(\ZZ, G)$ depends only on the mod-$\ell$ homotopy
type of $BG$. This appears as Corollary \ref{Z_hinv}. 

One would also like to know whether there is a result in the spirit of Friedlander--Mislin for $B(\ZZ,G)$. That is, let $G$ be
a Chevalley group, and consider the two associated groups $G(\bar{\FF}_p)$ and $G(\CC)$. The former is a locally
finite group, and the latter is a Lie group. We wish to compare the spaces $B(\ZZ, G(\bar\FF_p))$ and $B(\ZZ,
G(\CC))$. There is a zigzag of maps of spaces 
\begin{equation}
 \label{eq:6}
 \varphi:B(\ZZ,G(\bar{\FF}_p)) \leftrightarrow B(\ZZ,G(\CC)).
\end{equation}
We show that in general, this map is not a  mod-$\ell$ equivalence---here $\ell$ is a prime different from $p$. The failure
here is for substantive. One can produce $B(\ZZ, G)$ as a colimit over a topological poset of $BA$s, where the $A$s
range over some cofinal family of abelian subgroups of $G$. In saying this is a topological poset, we mean that the sets
of objects and of inclusions appearing carry a topology. For example, any Lie group $G$ has an associated space of
maximal tori, not merely a discrete set of such tori. In the colimit construction of $B(\ZZ, G)$, it is not sufficient
to consider the set of tori, the topology of the parametrizing space plays an important role. 

While a Friedlander--Milnor result holds for abelian groups, so that 
\[
\varphi: B A (\bar \FF_p) = B(\ZZ,A(\bar{\FF}_p)) \leftrightarrow B(\ZZ,A(\CC)) = B A( \CC)
\]
is a mod-$\ell$ equivalence, it is not possible to promote this to a mod-$\ell$ equivalence in \eqref{eq:6}, because on
the right hand side, the topology of the space of abelian subgroups intervenes in a significant way. On the left, the
collection of subgroups does not carry a comparable topology. In Example \ref{ex:noMF}, we show that $B(\ZZ,
\GL_2(\bar\FF_p))$ and $B(\ZZ, U(2))$ are not equivalent modulo $\ell$.

The topologies on spaces of abelian subgroups of compact Lie groups arise essentially as a consequence of the
continuous conjugation action of the group on itself. From this point of view, it is not perhaps wholly surprising that
when we take a Borel construction, i.e., we pass to homotopy orbits, for the conjugation action of $G$ on $B(\ZZ, G)$, the problem
caused by the innate topology on the poset of abelian subgroups disappears, and we are left with the content of Theorem
\ref{main_repeat}: there is a mod-$\ell$ equivalence in the spirit of \eqref{eq:6}, but only after taking homotopy
orbits. That is, there is a zigzag of mod-$\ell$ equivalences
 \begin{equation*}
 \varphi': EG(\bar{\FF}_p) \times_{G(\bar \FF_p)} (B(\ZZ,G(\bar{\FF}_p)) \leftrightarrow EG(\CC)\times_{G(\CC)}B(\ZZ,G(\CC)).
\end{equation*}
\section{Preliminaries} 
 
\subsection{Realization of simplicial spaces}
Throughout this paper, the term ``simplicial space'' will be taken to encompass ``simplicial set'' by endowing the sets
concerned with the discrete topology. Given a simplicial space $X$, we will often make use of the \textit{fat
 realization}, $\fat{X}$. This construction is similar to the usual realization, $|X|$, except that the degeneracy maps
are disregarded in the identification. Collapsing the degeneracies gives a natural map $\fat{X}\to |X|$. This map is a
homotopy equivalence when the simplicial space is \textit{good}, viz, when the degeneracy maps
$s_i:X_{n-1}\rightarrow X_{n}$ are closed cofibrations in the sense of Hurewicz, \cite[Proposition
A1]{SegalCategoriescohomologytheories1974a}. Fat realization is a homotopy colimit, and therefore a levelwise weak
equivalence $X\rightarrow Y$ induces a weak equivalence between the fat realizations.

\subsection{Mod-$\ell$ equivalences}

The topological spaces appearing in the sequel will have the homotopy types of CW complexes unless otherwise
noted---which is to say that they are not too poorly behaved. We will work throughout with spaces $X$
that are not nilpotent, that is, the structure of $\pi_1(X,x_0)$ and the action of this group on the higher
homotopy groups of the space may be complicated.

Fix a nonnegative integer $\ell$, generally a prime number. Write $\ZZ/\ell$ for the ring of integers modulo $\ell$.

\Def{\label{def:lequiv}\rm{We will say that a map $f: X \to Y$ of spaces is a \textit{mod-$\ell$ equivalence} if any of the following equivalent conditions is met
 \begin{enumerate}
 \item $f_* : H_*(X, \ZZ/\ell) \to H_*(Y, \ZZ/\ell)$ is an isomorphism;
 \item $f^*: H^*(Y, M) \to H^*(X, M)$ is an isomorphism for any $\ell$-complete abelian group $M$
 \item The Bousfield--Kan $\ell$-completion $f^{\wedge}_\ell : X^{\wedge}_\ell \to Y^{\wedge}_\ell$ is a weak equivalence.
 \end{enumerate}}}

The equivalence of these conditions follows from \cite[Lemma I.5.5]{BK72} and \cite[Chapter 10.1]{MP12}.

\Def{\rm{We will say that a map $f: X \to Y$ of spaces is a \textit{strong mod-$\ell$ equivalence} (see \cite{DZ87}) if the
 following conditions all hold:
 \begin{enumerate}
 \item $f_* : \pi_0(X) \to \pi_0(Y)$ is an isomorphism, and
 \item for any choice of basepoint $x \in X$, the map $f_*: \pi_1(X, x) \to \pi_1(Y, f(x))$ is an isomorphism, and
 \item for any choice of basepoint $x$, any of the evident maps $\tilde f: \tilde X_x \to \tilde Y_{f(x)}$ from the
 universal cover of the component of $x$ to the universal cover of the component of $f(x)$ is a mod-$\ell$ equivalence.
 \end{enumerate}}}

\subsection{Simplicial spaces of homomorphisms}
Let $\tau$ be a finitely generated discrete group. For a topological group $G$ the space of homomorphisms can be topologized as a subspace
\[
\Hom(\tau, G) \subset \prod_{g(\tau)} G
\]
where $g(\tau)$ denotes a set of generators. When $G$ is an affine algebraic group, $\Hom(\tau , G)$ has the structure of an
affine variety, \cite[Lemma 2.2]{Vil17}

Given a cosimplicial group $\tau^\dt$, the spaces of homomorphisms in
each simplicial degree fit together to produce a simplicial space
\[
\Hom(\tau^\dt,G),
\]
the structure maps being induced from the structure maps of the cosimplicial group. If $\tau^\dt$ is finitely
generated in each degree, we produce in this way a simplicial topological space. 

\Pro{\label{good} Let $G$ be an affine algebraic group defined over $k=\RR$ or $\CC$. 
 Then the simplicial space $\Hom(\tau^\dt,G(k))$ is good.}
\Proof{
The degeneracy maps are induced by surjective homomorphisms $s^i:\tau^n\rightarrow \tau^{n-1}$. By \cite[Proposition 1.7]{LM85} the induced map 
\[s_i:\Hom(\tau^{n-1},G)\rightarrow \Hom(\tau^n,G)\]
is a closed immersion. The triangulation theorem of \cite{Hir75} implies that the target space $ \Hom(\tau^n,G)$ has
the structure of a CW complex such that $s_i$ is the inclusion of a subcomplex. It follows that the degeneracy maps are
 closed cofibrations.}

\Def{\rm{ 
For a cosimplicial group $\tau^\dt$ and a topological group $G$ we define
\[ B(\tau,G) = |\Hom(\tau^\dt,G)|.\]
}}

 When $G$ is an affine algebraic group over $\CC$ (or $\RR$), up to homotopy we can use the fat realization as a consequence of Proposition \ref{good}.

It is also useful to work at a connected component. Let $\theta:\tau \rightarrow G$ denote a fixed homomorphism. The
component of the homomorphism space containing $\theta$ is denoted by $\Hom(\tau,G)_\theta$. In the simplicial context we can consider  $\theta:\tau^\dt\to G$, a homomorphism of cosimplicial groups where $G$ is given the trivial simplicial structure.  We can   define
$B(\tau,G)_\theta$ to be the geometric realization of the simplicial space $\Hom(\tau^\dt,G)_\theta$. Among all
components of $\Hom(\tau, G)$, there is a distinguished component: $\Hom(\tau, G)_1$, the component of the trivial homomorphism.

Let $\Delta$ denote the simplicial category, that is, the category having objects consisting of finite sets
$[n]=\{0, \dots, n\}$ and where the morphisms are order preserving functions. Let $\Delta^\dt$ denote the standard
cosimplicial simplicial set, where $\Delta^n$ is the usual $n$-simplex. There exists a \textit{standard reduced
 cosimplicial simplicial set} $\tilde \Delta^\dt$, where $\tilde \Delta^{n} = \Delta^n/ (\Delta^n_0)$, i.e., the
simplicial set obtained by identifying all vertices of $\Delta^n$ to a basepoint $\ast$. We define the \textit{reduced
 free cosimplicial group} $F^\dt$ to be $\pi_1(\tilde \Delta^\dt, \ast)$. For any abelian group $A$ we define the
\textit{reduced free abelian cosimplicial group with coefficients in $A$}, denoted $A^\dt$, by
$A^n = H_1( \tilde \Delta^n, A)$. The \textit{reduced free abelian cosimplicial group}, $\ZZ^\dt$, is used in the case
of $\ZZ$.

The definitions above have the advantage of being clearly functorial, but the disadvantage of being indirect. Since the
present paper relies heavily on these cosimplicial constructions, we illustrate how they work out in practice. Consider
the case of $A= \ZZ$, and fix a degree $n$. For each pair of integers $i,j$ satisfying $0 \le i \le j \le n$, we name
symbols $e_{i \to j}$, and let $\ZZ[n]$ denote the quotient of the free abelian group generated by symbols
$\{ e_{i \to j} \}_{i, j}$ subject to the relations $e_{i \to j} + e_{j \to k} = e_{i \to k}$. We remark that $e_{i \to i}$ is identified with $0$, but we nonetheless name $e_{i \to i}$ as a
technical convenience. There is a canonical isomorphism $\ZZ[n] \iso H_1( \tilde \Delta^n, \ZZ)$, the element
$e_{i \to j}$ mapping to the cycle corresponding to the edge $i \to j$. In the case where $i = j$, this edge is
degenerate and $e_{i \to i} = 0$. The relations on the cycles imposed by the higher cell structure of $|\tilde
\Delta^n|$ amount to $e_{i \to j} + e_{j \to k} = e_{i \to k}$. Therefore the $n$-th level of $\ZZ^\dt$ is canonically
isomorphic to $\ZZ[n] \iso \ZZ^n$.

It is possible to describe the cosimplicial structure of $\ZZ^\dt$ without difficulty in this language. Let $f:[n] \to
[m]$ be an order-preserving map, and consider $e_{i \to j} \in \ZZ[n]$, then $f_*(e_{i \to j}) = e_{f(i) \to f(j)}$. In
particular, fixing the ordered basis 
$\{e_{0\to 1}, e_{1 \to 2}, \dots, e_{n-1 \to n}\}$
 for each $\ZZ[n]$, the coface
maps are given by 
\[ \delta^i (a_1, \dots, a_n) = \begin{cases} (0,a_1, \dots, a_n) \text{ if $i=0$, }\\ (a_1, \dots, a_i, a_i, \dots, a_n)
 \text{ if $0 < i < n+1$}, \\ (a_1, \dots, a_n, 0) \text{ if $i=n+1$}. \end{cases} \]
For the codegeneracy maps:
\[ \sigma^i (a_1, \dots, a_n) = (a_1, \dots, \hat {a_i}, \dots a_n), \text{ where $\hat {a_i}$ denotes
 omission.} \]

 \medskip

 The above explicit constructions apply with the evident modification in the case of the cosimplicial group $(\ZZ/m)^\dt$. In the case of
 $F^\dt$, the same is also true with the caveat that the groups concerned are not abelian, and so care must be taken in
 the ordering of terms. Using the Hurewicz and Universal-Coefficient theorems, one obtains canonical surjections of
 cosimplicial groups
\begin{equation}\label{eq:sequence}
F^\dt \to \ZZ^\dt \to (\ZZ/m)^\dt. 
\end{equation} 

Given any cosimplicial group $\tau^\dt$ that is finitely generated in each degree and any topological group $G$, one may form
a simplicial space of maps $\Hom(\tau^\dt, G)$, which in degree $n$ is the space of homomorphisms $\Hom(\tau^n, G)$---since
the group $\tau^n$ is finitely generated, this may be topologized as a subspace of $G^{\times N}$ for some large
$N$. If $\tau^\dt$ happens to be the
reduced free cosimplicial group, $F^\dt$---without the abelian condition---then $\Hom(\tau^n, G)$ can be identified as
$G^{\times n}$, the set of $n$-tuples of elements in $G$, and the construction $B(\tau, G)$ recovers a familiar
construction of the classifying space $BG$ of the group $G$. When
$\ZZ^\dt$ is taken, $B(\ZZ, G)$ is the realization of the subobject of $BG_\dt$ obtained by restricting to $n$-tuples
of pairwise commuting elements in each degree, and $B(\ZZ/m, G)$ is the realization of the further-restricted
subobject where all group elements appearing must moreover be $m$-torsion. 
 The surjections in \ref{eq:sequence}
yield inclusions of spaces
\[ B(\ZZ/m,G) \rightarrow B(\ZZ,G) \rightarrow B(F,G)=BG.\]

\section{$B(\ZZ,G)$ as a homotopy colimit}

Let $B(\ZZ,G)_1$ denote the geometric realization of the simplicial space $\Hom(\ZZ^\dt,G)_1$ of path-connected components of homomorphism spaces at the trivial homomorphism $1:\ZZ^n\to G$. For compact connected Lie groups a homotopy colimit decomposition for $B(\ZZ,G)_1$ is given in \cite{AG15}. The goal of this section is to generalize this to a homotopy colimit decomposition for the full space $B(\ZZ,G)$ where $G$ is a compact Lie group, not necessarily connected.

\subsection{Homotopy colimits}

Let $\catA$ be a category internal to the category of topological spaces, i.e., a \textit{topological category} in the
language of \cite{Bor94}. We will denote the space of objects and space of morphisms by $\ob(\catA)$ and $\mor(\catA)$
respectively. The identity map $i:\ob(\catA)\rightarrow \mor(\catA)$, the source and the target maps
$s,t:\mor(\catA)\rightarrow\ob(\catA)$, and the composition map are all continuous. An \textit{internal functor} between two topological categories is 
a pair of continuous maps between the object and morphism spaces of the categories that respect the source, target, identity and composition maps.

There is also the concept of a space carrying a $\cat{A}$-action, called a $\catA$--\textit{module}. A left $\catA$--module is a space
$X$ together with a map $f_X:X\rightarrow \ob(\catA)$ and an action map $\mor(\catA)\times_{\ob(\catA) }X\rightarrow X$
that is unital and associative in the evident sense, \cite{HV92}. A \textit{morphism of left $\catA$--modules} is a
continuous map $X\rightarrow X'$ respecting the action of $\catA$.

There is also a notion of \textit{right $\catA$--module} defined in the obvious way. We usually denote a left $\catA$--module as a
functor $X:\catA\rightarrow \catTop$, and a right $\catA$--module as $Y:\catA^\op\rightarrow \catTop$. Under this correspondence the preimage $f_{X}^{-1}(a)$ is regarded as the value $X(a)$ of the functor at $a\in \ob(\catA)$, similarly for the contravariant functor $Y$. Given a left
$\catA$--module $X$ and a right $\catA$--module $Y$ we can define a simplicial space
\[ B(Y,\catA,X)_n = Y \times_{\ob(\catA)} \times \underbrace{\mor(\catA) \times_{\ob(\catA)} \cdots \times_{\ob(\catA)}
 \mor(\catA)}_n \times_{\ob(\catA)} X \]
with the usual face and degeneracy maps. The geometric realization of this space is called the \textit{bar
 construction} and will be denoted by $B(Y, \catA, X)$. The \textit{homotopy colimit} of the module $X$ is constructed as
$B(\ast,\catA,X)$ where $\ast:\catA \rightarrow \catTop$ is the constant functor. When $X$ is also the constant functor $\ast$ we write $B(\catA)$ instead of $B(\ast,\catA,\ast)$. 
 
\Rem{\label{cofibration}\rm{ We shall assume that the topological categories we consider satisfy the property that the
 identity map $i:\ob(\catA) \rightarrow \mor(\catA)$ is a closed cofibration over $\ob(\catA)\times \ob(\catA)$ with
 respect to the diagonal map and the source-target map $s\times t$. This guarantees that the simplicial space
 $B(Y,\catA,X)$ is proper (see \cite[Rem. 2.3]{Mey86} or \cite[Def. A.1]{Lin13}).}}

Given an internal functor $F:\catA'\to \catA$ the undercategory $a \downarrow F$, i.e., the category in which the
objects are pairs $x \in \ob(\catA')$ and $f: a \to F(x) \in \mor(\catA)$ is a topological category in an evident
way: the space of objects is a subspace of $\ob(\catA')\times \mor(\catA)$, and the space of morphisms inherits a
topology from $\mor(\catA')$.  We say $F$ is \textit{cofinal} if $B(a\downarrow F)$ is contractible for all
$a\in \ob (\catA)$.  We will refer to the next result quite often.

\Pro{\cite[Lemma A.5]{Lin13} \label{cofinality}
Let $X$ be a $\catA$-module and $F:\catA'\to \catA$ be a cofinal functor. Then the map
$$
\hocolim{\catA'} X\circ F \to \hocolim{\catA} X
$$
is a homotopy equivalence.
} 

We introduce some notation as a preparation for our main result concerning homotopy colimits. We remind the reader that
$\catDelta$ denotes the category having as objects the finite sets $[n] = \{0, 1, \dots, n\}$ and where the morphisms
are order-preserving maps. The category $\cat{Cat}$ is the category of small categories. There is a functor $\Delta^\dt: \catDelta \rightarrow \Cat$ sending $[n]$ to the category $0\rightarrow
1\rightarrow \cdots \rightarrow n$.

\Def{ \rm{ \label{DefOverCat}
Given a small category $\catC$, the overcategory $\Delta^\dt\downarrow \catC$ is described as follows:
\begin{itemize}
\item The objects of $\Delta^\dt\downarrow \catC$ consist of pairs $([k], \sigma: [k] \to \catC)$.
\item The morphisms $([k_0], \sigma_0) \to ([k_1], \sigma_1)$ are order-preserving maps $\theta: [k_1] \to [k_0]$ such
 that $\sigma_1 = \sigma_0 \theta$.
\end{itemize}
There is a forgetful functor $U:\Delta^\dt\downarrow\catC \rightarrow \catDelta^\op$.
}}

Recall that $\catA$ is a topological category. Let $X$ be a left $\catA$--module. Let $\catC$ be a small category, regarded as
a topological category by giving the sets of objects and morphisms the discrete topology. Let $\rho: \catA \to \catC$ be
an internal functor. We may define a functor
\[ \rho_*X:\Delta^\dt\downarrow\catC \to \catTop \]
as follows.

There is a composite map 
\[ \tilde X_k : B(\ast, \catA, X)_k \to B(\ast, \catA,\ast)_k \overset{B \rho}{\to} B(\ast,\catC,\ast)_k, \]
here $B(\ast, \catC, \ast)_k$ is the set of functors $\sigma: [k] \to \catC$, and therefore any such functor $\sigma$
has a preimage $\tilde X_k^{-1}(\sigma) \subset B(\ast, \catA, X)_k$. We define 
\begin{equation}
 \label{eq:8}
 \rho_* X(\sigma) = \tilde X_k^{-1}(\sigma).
\end{equation}
Explicitly, if $\sigma$ is the data of a sequence $c_0 \to c_1 \to \dots \to c_k$ of maps in $\catC$, then
$$
\rho_*X(\sigma)=\set{(a_0\to\cdots\to a_k,x)|\; \rho(a_i)=c_i\;\text{ and }\; x\in X(a_0)}
$$
with the subspace topology induced from $B(\ast, \catA, X)_k$.
Given a morphism $([k_0], \sigma_0) \to ([k_1], \sigma_1)$ induced by an ordinal map $\theta: [k_1]\rightarrow [k_0]$ the functor $\rho_*X$ sends the pair $(a_0 \to a_1 \dots \to a_{k_0},x)$ where $x\in X(a_0)$ to the pair $(a_{\theta(0)} \to a_{\theta(1)} \dots \to a_{\theta(k_1)},X(\alpha)(x))$ where $X(\alpha):X(a_0)\to X(a_{\theta(0)})$ is induced by the composite map $\alpha:a_0\rightarrow a_{\theta(0)}$.

Our main observation is the following result.

\begin{thm} \label{HocolimTopToDisc}
Let $\catA$ be a topological category and $\catC$ be a discrete category. 
Suppose we are given a diagram 
$$
\begin{tikzcd}
\catA \arrow{d}{\rho} \arrow{r}{X} & \catTop \\
\catC
\end{tikzcd}
$$
consisting of an $\catA$-module $X$ and an internal functor $\rho$. Then there is a natural homotopy equivalence
$$
\hocolim{\Delta^\dt\downarrow \catC} \rho_*X \to \hocolim{\catA} X.
$$
\end{thm}

The virtue of this theorem is that the homotopy colimit over a topological category is replaced by a homotopy colimit over a discrete category.

The proof is broken into two lemmas.

\Lem{\label{proper} There is a natural homotopy equivalence
\[
B(\ast, \Delta^\op, B(\ast,\catA, X)_\dt) \rightarrow B(\ast, \catA,X).
\]
}
\Proof{This map is given by the Bousfield--Kan map $B(\ast,\Delta^\op,-)\rightarrow |-|$.
Under the cofibrancy assumption on $\catA$ given in Remark \ref{cofibration} the bar construction $B(\ast,\catA, X)_\dt$
is a proper simplicial space. Recall from \cite{Stromhomotopycategoryhomotopy1972} that there is a \textit{Str\o m model structure} on the category of
topological spaces, for which the fibrations are the Hurewicz fibrations, the cofibrations are the closed Hurewicz
cofibrations and where the weak equivalences are homotopy equivalences in the usual sense. 

When we use the Str\o m model structure on the category of topological spaces the Reedy cofibrant simplicial objects coincide with proper simplicial spaces. The Bousfield--Kan map is a natural weak equivalence for Reedy cofibrant simplicial objects \cite[Theorem 18.7.4]{Hir03}.
}

A functor $\sigma:[n]\rightarrow\Delta^\dt\downarrow\catC$ consists of a sequence $\sigma_0 \to \cdots \to\sigma_n$ of functors $\sigma_i: [k_i]\to \catC$ together with ordinal maps $\theta_i:[k_{i+1}]\to [k_{i}]$ such that $ \sigma_{i+1}=\sigma_i\theta_i$. Observe that $\sigma$ is determined by $\sigma_0$ and the ordinal maps $\theta_i$ for $0\leq i\leq n-1$.

For $n\geq 0$ we define a map
\[
\beta_n:\coprod_{\sigma:[n]\rightarrow\Delta^\dt\downarrow\catC} \rho_*X(\sigma_0) \rightarrow \coprod_{\sigma':[n]\rightarrow \catDelta^\op} B(\ast,\catA,X)_{\sigma'(0)}
\]
induced by the inclusion $\rho_*X(\sigma_0)\subset B(\ast,\catA,X)_{k_0}$ at the term corresponding to $\sigma$. Let $\beta$ denote the map $|\beta_\dt|$.

\Lem{\label{quotient} 
There is a natural homeomorphism 
\[\beta:B(\ast,\Delta^\dt\downarrow\catC,\rho_*X)\rightarrow B(\ast, \catDelta^\op,B(\ast,\catA,X)_\dt) .\]
}
\Proof{ 
For a fixed $\sigma'$ given by the sequence of ordinal maps $[l_0]\leftarrow [l_1] \leftarrow \cdots \leftarrow [l_n]$ consider the preimage of $B(\ast,\catA,X)_{l_0}$ under $\beta_n$.
Since $\sigma$ is determined by $\sigma_0$ and the sequence of ordinal maps $\theta_i$ the restriction of $\beta_n$ to the preimage is given by
\[
\coprod_{\gamma: [l_0]\rightarrow \catC} \tilde{X}^{-1}_{l_0}(\gamma) \rightarrow
B(\ast,\catA,X)_{l_0}
\]
 which is a homeomorphism by the definition of $\tilde{X}_{l_0}$.
}

\begin{example}\label{ex:hocolim}
 A special case of Theorem~\ref{HocolimTopToDisc} that occurs often in the sequel is the following simple
 situation. Suppose we have two spaces $A_0$ and $A_1$ of objects, and a space of nontrivial morphisms $M$ from $A_0$
 to $A_1$, i.e., the set $M$ is equipped with a source map $s: M \to A_0$ and a target map $t: M \to A_1$. We package
 this all as a topological category $\cat A$ by setting $A_0 \amalg A_1$ to be the space of objects and $M \amalg A_0
 \amalg A_1$ the set of morphisms, with $A_0 \amalg A_1$ comprising the identity maps.

 This topological category is a topologized version of the discrete category $[1] = (0 \to 1)$, and is equipped with an
 obvious internal functor $\rho : \cat A \to [1]$, sending $A_i \mapsto i$ and $M$ to the non-identity morphism.

 An $\cat A$-module $X$ amounts to two spaces $X_0$ and $X_1$, each equipped with a map $X_i \to A_i$, and such that there
 is a map $M \times_{A_0} X_0 \to X_1$.

 We might wish to calculate $\hocolim{\cat A} X$, and to do so, we use Theorem \ref{HocolimTopToDisc}. According to
 the theorem, we should calculate a homotopy colimit of an infinite diagram, the image of a functor on
 $\Delta^\bullet \downarrow [1]$, but it will turn out that almost all terms in this diagram are degenerate and can be
 disregarded.

 The objects of the category $\Delta^\bullet \downarrow [1]$ are the functors $[k] \to [1]$, and the integer $k$ gives an
 obvious grading on these objects. In degree $0$ there are two objects, corresponding to the two functors $f_0, f_1 : [0] \to
 [1]$, the first with image $0$ and the second with image $1$. To these we associate $\rho_*(X)(f_0) = X_0$ and
 $\rho_*(X)(f_1) = X_1$.

 In degree $1$, there are three objects: the identity map $\mathrm{id}_{[1]}: [1] \to [1]$ and two degenerate maps $s_0, s_1: [1] \to
 [1]$. To the nondegenerate map, we associate $\rho_*(X)(\mathrm{id}_{[1]}) = M \times_{A_0} X_0$. To the degenerate map $s_0$
 at $0$ we associate $A_0 \times_{A_0} X_0 \approx X_0$ and to the degenerate map $s_1$ at $1$ we associate $A_1 \times_{A_1}
 X_1 \approx X_1$---in each case $A_i$ is serving both as the space of objects and of morphisms. Up to degree $1$,
 then, the diagram $\rho_*(X): \Delta^{\le 1} \downarrow [1] \to \cat{Top}$ is given by 
 \[ \xymatrix{ X_0 \ar^\iso@<0.7ex>[dr] \ar_\iso@<-0.7ex>[dr] & & M \times_{A_0} X_0 \ar[dr] \ar[dl] && X_1 \ar^\iso@<0.7ex>[dl] \ar_\iso@<-0.7ex>[dl] \\ & X_0 & & X_1 } \]
 The maps $X_i \to X_i$ are identity maps, the map $M \times_{A_0} X_0 \to X_0$ is projection and the map $M
 \times_{A_0} X_0 \to X_1$ is the map demanded by the $\cat A$-module structure on $X$.

 In higher degrees, $d \ge 2$, all objects appearing are degenerate, which is to say, of the form $(a_0 \to \dots \to a_k, x)$
 where at least one arrow appearing is an identity arrow. In the diagram, all maps emanating from degenerate objects
 are identity maps, and these objects may therefore safely be disregarded when calculating the homotopy colimit. In
 summary, the original topological homotopy colimit is equivalent to the discrete homotopy colimit of the diagram
 \begin{equation*}
 \xymatrix{ M \times_{A_0} X_0 \ar[r] \ar[d] & X_1 \\ X_0 }
 \end{equation*}

\end{example}

We remark that it is not in principle difficult to generalize this example to the case of a topological directed
category, e.g., one modelled on $0 \to 1 \to \dots \to n$ rather than simply on $0 \to 1$, but the extra bookkeeping
required means that the general case does not serve so well as an example.

\subsection{Poset of abelian subgroups} Let $G$ be a compact Lie group. Let $\aA(G)$ denote the poset of all closed
abelian subgroups of $G$, partially ordered by inclusion. We topologize $\aA(G)$ by the Chabauty topology. In general,
this is one of two obvious choices for a topology on $\aA(G)$, the other being the finite topology. For compact groups
the topologies coincide, see \cite{BHK09}. In the Chabauty topology, $\aA(G)$ is a compact Hausdorff space for which
conjugation is continuous. We work more generally with a subcollection $\aA\subset \aA(G)$ closed under conjugation and
taking intersections. We say $\aA$ is \textit{cofinal} in $\aA(G)$ if the inclusion functor $\aA \to \aA(G)$ is a cofinal
internal functor.

Let $\bar{\aA} $ denote the set of conjugacy classes of subgroups in $\aA$ under the conjugation action by $G$.
We give $\bar \aA$ the quotient topology induced from $\aA$. We regard $\bar \aA$ as a poset and write $[A]\to [B]$ whenever $A$ is conjugate to a subgroup of $B$. Identifying conjugate subgroups gives an internal functor 
$$\rho:\aA \to \bar \aA.$$

Let $\tau^\dt$ be a cosimplicial group that is
a finitely generated abelian group in each degree. We define a left $\aA$--module $B(\tau, -): \aA \to \catTop$ by
defining the space of $B(\tau, -)$ as
\[ \bigcup_{A \in \ob(\aA)} B(\tau, A) \]
given the subspace topology as a subspace of $\ob(\aA) \times BG$. The action of $\mor(\aA)$ on $B(\tau, -)$ is on
the first coordinate.

Associated to this module there is a functor
$\rho_*B(\tau,-): \Delta^\dt\downarrow\bar\aA \rightarrow \catTop$, which sends $\bar\sigma:[k]\rightarrow\bar\aA$ to
the preimage under the map
\[
B(\ast,\aA,B(\tau,-))_{\dt} \rightarrow B(\bar\aA)_{\dt} 
\]
as defined in \eqref{eq:8}, where the role of $X$ is played by $B(\tau,-)$.

The following theorem allows one to calculate $B(\ZZ, G)$ by means of a homotopy colimit over
a discrete category. The idea is essentially due to \cite[Section 5]{AG15} where $G$ is assumed to be connected.

\Pro{\label{decompositionCofin} Let $G$ be a compact Lie group and $\aA$ be a subcollection of $\aA(G)$ closed under conjugation and taking intersections.
Suppose that the topological poset $\aA$ is cofinal in $\aA(G)$ and the quotient $\bar\aA$ is a discrete category. Then the map 
\[
\hocolim{\Delta^\dt\downarrow\bar\aA} \rho_*B(\tau,-) \rightarrow B(\tau,G)
\]
natural in both $\tau^\dt$ and $G$ is a weak equivalence.
}
\Proof{
We apply Theorem \ref{HocolimTopToDisc} with $X=B(\tau,-)$ to obtain a homotopy equivalence
$$
\hocolim{\Delta^\dt\downarrow\bar\aA} \rho_*B(\tau,-) \to \hocolim{\aA} B(\tau,-).
$$
It suffices to show that
\begin{equation}\label{mapWeq}
B(\ast,\aA,B(\tau,-)) \rightarrow B(\tau,G)
\end{equation}
is a weak equivalence. Let $\catA$ denote the topological category given by the Grothendieck construction of
$B(\tau,-)$, and let $\catA'$ be the topological category with object space $B(\tau,G)$ and only identity
morphisms. The object space of $\catA$ consists of pairs $(A,x)$ where $A\in \aA$ and $x\in B(\tau,A)$ topologized as a subspace of $\ob(\aA) \times BG$. Morphisms $(A,x)\to (A',x')$ are given by morphisms $\phi:A\to A'$ in $\aA$ such that $\phi(x)=x'$. This is a natural subspace of $\mor(\aA)\times BG$.
There is an internal functor $F:\catA\rightarrow \catA'$ such that $BF$ can be identified with the map in
(\ref{mapWeq}). Proposition \ref{cofinality} finishes the proof provided we can show that $F$ is cofinal.

 The undercategory
$x\downarrow F$ for an element $x \in B(\tau,G)$ consists of pairs $(A,x)$ such that $B(\tau,A)$ contains $x$. Since
$\aA$ is closed under intersections, the category $x \downarrow F$ has an initial element. In the topological context the map
$\ob(x\downarrow F)\rightarrow \mor(x \downarrow F)$ which sends an object to the unique morphism from the initial
object is required to be continuous, \cite{Lin13}. This condition is automatically satisfied for topological posets. Thus $B(x\downarrow F)$ is contractible.
}

Next we show how to construct a subcollection with the desired properties.
Fix an embedding $\iota:G\rightarrow U(n)$ into a unitary
group. Let $\tT$ denote the set of all intersections of maximal tori in $U(n)$. We define a collection of subgroups
\[
\aA_\iota =\set{G\cap S|\; S\in \tT} \subset \aA(G).
\] 
Let $\bar{\aA}_\iota$ denote the set of conjugacy classes of subgroups in $\aA_\iota$ under the conjugation action by
$G$. 
We will usually drop the embedding from the notation and simply denote these posets by $\aA$ and $\bar\aA$.

The conjugation action of $G$ on $\aA$ extends in the obvious way to functors $\sigma: [k] \to \aA$. 
This action is continuous with respect to the topology on $B(\aA)_k$. Let $N(\sigma)$
denote the stabilizer of $\sigma$ under the $G$-action.

\Lem{\label{QuoFiniteSet} The quotient space $B(\aA)_k/G$ is a finite set with discrete topology. In particular, $\bar\aA$ is finite.}
\begin{proof} Let $d_0:B\aA_k\to B\aA_{k-1}$ denote the map which forgets the first subgroup in a simplex.
Consider the pullback 
$$
\begin{tikzcd}
\aA\downarrow A_1 \arrow{r} \arrow{d} & B\aA_k \arrow{d}{d_0} \\
\set{\sigma} \arrow{r} & B\aA_{k-1}
\end{tikzcd}
$$
along a simplex $\sigma$ given by $A_1\to\cdots\to A_k$. First we will show that $\bar\aA$ is finite ($k=0$ case) and do induction on $k$. Then from the diagram we see that it suffices to show that $(\aA\downarrow A_1)/N(\sigma)$ is finite since by the induction step $(B\aA_{k-1})/G$ is finite.
Both of these claims follow from the general fact: Let $H\subset G$ be a closed subgroup. If a collection of subgroups of $H$ falls into finitely many $G$-conjugacy classes then the set of $H$-conjugacy classes is also finite, see \cite[Lemma 6.3]{Q71}. For the base case we apply this to the inclusion $G\subset U(n)$ and use the fact that $\tT/U(n)$ is finite (\cite[Theorem 5.4]{AG15}). Next we consider the inclusion $N(\sigma)\subset G$ and the collection $\aA\downarrow A_1$.
We conclude that $(B\aA_k) /G$ is finite as a set.
Note that being a quotient of a Hausdorff space by a compact group action this set has discrete topology.
\end{proof}

We identify 
\[
\aA = \coprod_{[A]\in \bar \aA} G/NA.
\]
Since the right hand side is a topological space, we may declare $\aA$ to be a topological poset with this space of
objects and where the space of morphisms consists of the subspace of $\aA \times \aA$ of pairs $(A_1, A_2)$ where $A_1
\subseteq A_2$. We observe that the topology on $\aA$ is compact and Hausdorff. Moreover, this topology coincides with the subspace topology induced from $\aA(G)$.

\Rem{\rm{
In particular, the identity map $i:\ob(\aA) \rightarrow \mor(\aA)$ satisfies the cofibration condition in Remark \ref{cofibration}.
}}

Next we give an explicit description of the functor $\rho_*B(\tau,-)$ with respect to a fixed embedding $\iota:G\to U(n)$.

The reduction functor $\rho$ induces a functor between the $k$--simplices
\[
B(\rho)_k: B(\aA)_k \rightarrow B(\bar\aA)_k
\]
of the classifying spaces. Let $Q(\bar\sigma)$ denote the preimage $B(\rho)_k^{-1}(\bar\sigma)$ of a $k$--simplex
$\bar{\sigma}$. 

The functor $\rho_*B(\tau,-):\Delta^\dt \downarrow \bar\aA \to \catTop$ sends a simplex $\bar\sigma$ to the space
$$
\rho_*B(\tau,\bar\sigma) = \bigcup_{\sigma\in Q(\bar\sigma)} B(\tau,\sigma(0))
$$
(see \ref{eq:8}). Let $\sigma$ be a simplex in the preimage $Q(\bar\sigma)$. There is a map
\begin{equation}\label{mapConj} 
G\times B(\tau,\sigma(0)) \to \rho_*B(\tau,\bar\sigma)
\end{equation}
defined by conjugation $(g,x)\mapsto gxg^{-1}$. 
 The stabilizer group $N(\sigma)$ acts on $G\times B(\tau,\sigma(0))$ by $n(g,x)=(gn^{-1},nxn^{-1})$. The conjugation map (\ref{mapConj}) is equivariant with respect to this action, and it factors through the quotient by the action of $N(\sigma)$.

\Lem{\label{homeo} The conjugation map
\[
\coprod_{[\sigma]\in Q(\bar\sigma)/G} G\times_{N(\sigma)} B(\tau,\sigma(0)) \rightarrow \rho_*B(\tau,\bar\sigma)
\]
is a homeomorphism for all $\bar\sigma$, and $Q(\bar\sigma)/G$ is a finite set.
}
\begin{proof}
As a consequence of Lemma \ref{QuoFiniteSet} the quotient space $Q( \bar\sigma)/G$ is a finite set. Note that there is a continuous surjective map $q:\rho_*B(\tau,\bar\sigma)\to Q(\bar\sigma)$. We can compose this with the quotient map $Q(\bar\sigma)\to Q(\bar\sigma)/G$ to obtain the map $\bar q: \rho_*B(\tau,\bar\sigma)\to Q(\bar\sigma)/G$. Therefore it suffices to show that each factor of the disjoint union maps homeomorphically onto the corresponding component of $\rho_*B(\tau,\bar\sigma)$ over a point $[ \sigma]$. In details, the restriction of $q$ to this component gives us an equivariant map $\bar q^{-1}[\sigma] \to G/N(\sigma)$.
Now the result follows from the general fact that for a $G$-space $X$ together with a $G$-map $f:X\to G/H$ the natural map $G\times_H f^{-1}(H)\to X$ is a homeomorphism.
\end{proof}

Combining Proposition \ref{decompositionCofin} with Lemma \ref{QuoFiniteSet} we obtain the main result of this section. 
 
\Thm{\label{decomposition}
Let $G$ be a compact Lie group and let $\tau$ be a cosimplicial group such that $\tau^n$ is a finitely generated abelian
group for all $n\geq 0$. Fix an embedding $\iota: G \to U(n)$ and let $\bar \aA$ denote the set of $G$ conjugacy classes
of closed abelian subgroups given by intersecting $G$ with intersections of maximal tori of $U(n)$. The natural map 
\[
\hocolim{\Delta^\dt\downarrow\bar\aA} \rho_*B(\tau,-) \rightarrow B(\tau,G)
\]
is a weak equivalence.
}

\Rem{\rm{\label{HomSpaceDec} 
Proposition \ref{decompositionCofin} also holds for the functor $\Hom(\tau^n,-):\aA \to \catTop$. Thus we have a similar homotopy colimit decomposition for the homomorphism spaces which gives us a weak equivalence
\[
\hocolim{\Delta^\dt\downarrow\bar\aA} \rho_*\Hom(\tau^n,-) \rightarrow \Hom(\tau^n,G).
\]
}}

\section{Mod-$\ell$ homology of $B(\ZZ,G)$}
Fix a prime $\ell$.
Suppose $A$ is a compact abelian Lie group and let $A_{\ell}$ denote the $\ell$-primary torsion subgroup of $A$. It is
the case that $BA_{\ell} \to BA$ is a mod-$\ell$ equivalence. The proof of this fact is elementary when $A$ is finite,
and follows from the Milnor conjecture for $U(1)$. The general case follows from these two cases
without difficulty.
In this section, we generalize this observation to two classes of groups: compact Lie groups and locally finite groups. That is, if $G$ is a group of either of these two classes, then
we establish that the natural map
\begin{equation}
 \label{eq:1}
 \colim{k \to \infty} B( \ZZ/\ell^k , G) \to B(\ZZ, G) 
\end{equation}
is a mod-$\ell$ equivalence.
The case of finite $G$ already appears in \cite[Theorem 3.4]{O14} 

As an example we consider the double cover map $R:\SU(2)\to \SO(3)$. By explicit calculation we show that $R$ induces a mod-$\ell$ equivalence
$$
B(\ZZ,\SU(2)) \to B(\ZZ,\SO(3))
$$ 
for all $\ell>2$, which is expected as a consequence of the $\ell^k$-torsion abelian subgroup structures of the underlying compact Lie groups.
 
\subsection{Mod-$\ell$ homology} We will use the homotopy colimit description given in Theorem \ref{decomposition} to describe the mod-$\ell$ homology of $B(\ZZ,G)$. First we rewrite the colimit that appears in \eqref{eq:1} in a more appealing way. 

Let $\ZZ_\ell$ denote the additive group of $\ell$-adic integers. This group has the profinite topology. Given a Lie group $G$ we can consider the set of continuous homomorphisms $\Hom(\ZZ_\ell,G)$ in the category of topological groups. 

\Lem{\label{Hom_ell_adic}
Let $G$ be a Lie group or a discrete group. Then there is a natural bijection
$$
\colim{k} \Hom((\ZZ/\ell^k)^n,G) \to \Hom((\ZZ_\ell)^n,G)
$$
for all $n\geq 0$. } 
\Proof{This statement holds for any profinite group, and follows from the fact that Lie groups do
not contain ``small subgroups". There is a neighbourhood $U$ of the identity of $G$ so that the only subgroup contained
in it is the trivial subgroup. It suffices to focus on $n=1$.  Let $f:\ZZ_\ell\to G$ be a continuous homomorphism. 
The preimage $f^{-1}(U)$ is an open set and thus contains an open subgroup $V$. By the assumption on $U$ the subgroup $V$ is contained in the kernel of $f$. This implies that the kernel is open. By compactness of the
profinite group the image of $f$ is finite. Therefore $f$ factors through one of the quotients $\ZZ/\ell^k$.

The case when $G$ is discrete also follows easily from the compactness of $\ZZ_\ell$
}

\Def{\rm{
Let $\Hom((\ZZ_\ell)^n,G)$ denote the set of continuous homomorphisms with the topology induced from the colimit of the spaces $\Hom((\ZZ/\ell^k)^n,G)$ over $k$. We denote by $B(\ZZ_\ell,G)$ the geometric realization of the simplicial space $\Hom(\ZZ_\ell^\dt,G)$.
}
}

In effect we identify
\[
B(\ZZ_\ell,G)=\colim{k\to \infty} B(\ZZ/\ell^k,G).
\]
The completion map $\ZZ\to\ZZ_\ell$ at the prime $\ell$ induces a natural map $B(\ZZ_\ell,G)\to B(\ZZ,G)$.

\Cor{\label{lie_modulo}Let $G$ be a compact Lie group. The natural map
\[\alpha: B(\ZZ_\ell,G) \rightarrow B(\ZZ,G)\]
is a mod-$\ell$ equivalence.
}
\Proof{
Up to homotopy we can use the fat realization in the definition of $B(\ZZ/\ell^k,G)$. Then the maps in the direct limit are cofibrations, hence we can replace the direct limit by a homotopy direct limit. Using Theorem \ref{decomposition} and commuting the homotopy colimits
\[
\hocolim{k} \hocolim{\Delta^\dt\downarrow \bar\aA} \rho_*B(\ZZ/\ell^k,-) \simeq \hocolim{\Delta^\dt\downarrow \bar\aA}\hocolim{k} \rho_*B(\ZZ/\ell^k,-) 
\]
we are reduced to comparing 
\[
\hocolim{k} \rho_*B(\ZZ/\ell^k,\bar\sigma) \rightarrow \rho_*B(\ZZ,\bar \sigma)
\]
for an object $\bar\sigma$ in $\Delta^\dt\downarrow\bar\aA$.
The map in Lemma \ref{homeo} is natural with respect to $\tau$. Thus, it suffices to consider the map
\[
 G\times_{N(\sigma)} B(\ZZ/\ell^k, \sigma(0)) \rightarrow G\times_{N(\sigma)} B(\ZZ, \sigma(0))
\]
which is a map of fibrations over $G/N(\sigma)$ with the inclusion map
$$B(\ZZ/\ell^k, \sigma(0))\rightarrow B(\ZZ, \sigma(0))$$ between the fibres. This reduces the problem to the case of
the product of a torus with a finite abelian group. We can even restrict to each case separately since both $B(\ZZ,-)$
and $B(\ZZ/\ell^k,-)$ commutes with products. For finite abelian groups this is easy to verify since mod-$\ell$
cohomology is determined by the $\ell$-torsion part of the group. For tori it suffices to consider the rank one
case. When $G=U(1)$ the space $B(\ZZ/\ell^k,U(1))$ can be identified with $B\ZZ/\ell^k$, and the map $\alpha$ can be
identified with the map
\[
B(\ZZ/\ell^\infty) \rightarrow B U(1)
\]
induced by $\ZZ/\ell^k\cong \mu_{\ell^k} \rightarrow U(1)$. The induced map in mod-$\ell$ cohomology is an isomorphism as a consequence of 
\[
H^*(B\ZZ/\ell^\infty ,\ZZ/\ell) = \ilim{k}{} H^*(B\ZZ/\ell^k ,\ZZ/\ell) = \ZZ/\ell[t]
\]
where $t$ is a degree $2$ generator. 
}

\subsection{Locally finite groups} 
Corollary \ref{lie_modulo} is an extension of the finite group case proved in \cite[Theorem 3.4]{O14}. We show that this easily generalizes to locally finite groups, that is, groups of the form 
\[
H=\colim{i} H_i
\]
where each $H_i$ is a finite group. We specifically will apply this in the case of the locally finite groups
$G(\bar{{\FF}}_p)$, the colimit of the finite groups $G(\FF_{p^i})$.

\Pro{\label{loc_finite}Let $H$ be a locally finite group. The natural map
\[
 B(\ZZ_\ell,H) \rightarrow B(\ZZ,H)
\]
is a mod-$\ell$ equivalence.
}
\Proof{ $B(\ZZ_\ell,G)$ is identified with the colimit of $B(\ZZ/\ell^k,G)$ over $k$ (Lemma \ref{Hom_ell_adic}).
 Note that for any finitely generated group $\pi$ the natural map
\[
\colim{i} \Hom(\pi,H_i) \rightarrow \Hom(\pi,H)
\]
is a bijection. Moreover, the geometric realization functor commutes with colimits. Therefore we have a sequence of homeomorphisms
\[
\colim{k} B(\ZZ/\ell^k,H) \approx \colim{k} \colim{i} B(\ZZ/\ell^k,H_i) \approx \colim{i} \colim{k} B(\ZZ/\ell^k,H_i)
\]
where the last one is a consequence of the commutativity of colimits.
By naturality the result follows from the fact that
\[
\colim{k} B(\ZZ/\ell^k,H_i) \rightarrow B(\ZZ,H_i)
\]
induces a mod-$\ell$ equivalence.
}

\subsection{Calculations}  
We make some preliminary observations to make the homotopy colimit decomposition in Theorem \ref{decomposition} useful
in practice. The cofinality result, Proposition \ref{cofinality}, allows us to work over simpler categories. The
collection $\aA$ that we use to decompose $B(\ZZ,G)$ can be replaced by a smaller collection. Let $\aA'\subset \aA$
denote the subcollection of subgroups containing the centre $Z\subset G$ and give $\aA'$ the subspace topology. The
inclusion functor $F:\aA'\to \aA$ is an internal functor. For any object $A$ of $\aA$ the classifying space of the undercategory
$A\downarrow F$ is contractible since the subgroup $AZ$ is initial in the category. Thus $F$ is cofinal.

Another simplification is to replace the simplex category $\Delta^\dt\downarrow \bar\aA$ with the full subcategory of non-degenerate simplices as already discussed in a special case in Example \ref{ex:hocolim}. Let $(\Delta^\dt\downarrow \bar\aA)_{\text{non-deg}}$
denote the category whose objects are non-degenerate simplices of $\bar \aA$, and morphisms are inclusions of simplices. For a simplex $\bar\sigma$ there is a unique non-degenerate simplex $\bar\sigma_0$ with a morphism $\bar\sigma_0 \to \bar\sigma$. This assignment is natural and defines a functor $R:\Delta^\dt\downarrow \bar\aA\to (\Delta^\dt\downarrow \bar\aA)_{\text{non-deg}}$, the right adjoint of the inclusion functor $L:(\Delta^\dt\downarrow \bar\aA)_{\text{non-deg}} \to \Delta^\dt\downarrow \bar\aA$ (under the conventions of Definition \ref{DefOverCat}). By the description of the functor $\rho_*B(\tau,-)$ in Lemma \ref{homeo} we see that its values at $\bar\sigma_0$ and $\bar\sigma$ are the same. 
Therefore we have $\rho_*B(\tau,-)= \rho_*B(\tau,-) L R $ and we observe that $R$ is cofinal.
 
In our examples $(\Delta^\dt\downarrow \bar\aA)_{\text{non-deg}}$ will be related to the poset $d[n]$ of non-empty subsets of $[n]$ ordered by inclusion. 
 
\subsection{$\SU(2)$} We think of $\SU(2)$ as the group of unit quaternions in the quaternion algebra $\mathbb{H}$ with the basis $\set{\one,\ii,\jj,\kk}$ where
$$
\one=\begin{bmatrix}
1&0\\
0&1
\end{bmatrix}\;\;
\ii=\begin{bmatrix}
i&0\\
0&-i
\end{bmatrix}\;\;
\jj=\begin{bmatrix}
0&i\\
i&0
\end{bmatrix}
$$ 
and $\kk=\ii\jj$. Let $T$ denote the maximal torus consisting of diagonal matrices, which is isomorphic to $\U(1)$. The normalizer $N(T)$ is generated by $T$ and $\jj$. Two different maximal tori intersect at the centre $Z=\set{\pm \one}$. 
The object space of $\aA$ is homeomorphic to $\SU(2)/N(T)_+$.
The base point corresponds to the centre, and the rest parametrizes the maximal tori. The category $\bar\aA$ can be identified with $[1]=\set{0\to 1}$, and $(\Delta^\dt\downarrow\bar\aA)_{\text{non-deg.}}$ with the opposite category $d[1]^\op$. As in \cite{AG15} we obtain the homotopy push-out diagram
\[
\begin{tikzcd}
\SU(2)\times_{N(T)}BZ \arrow{r}\arrow{d} & \SU(2)\times_{N(T)}BT \arrow{d}\\
BZ \arrow{r} & B(\ZZ,\SU(2))
\end{tikzcd}
\]
\subsection{$\SO(3)$}
The case of $\SO(3)$ is more complicated as there are abelian subgroups not conjugate to a maximal torus. Let
$R:\SU(2)\to\SO(3)$ denote the covering homomorphism.

The image $\bar T =R(T)$ is a maximal torus isomorphic to $\SO(2)$. Other than the trivial subgroup there are also finite abelian subgroups not conjugate to $\bar T$. But all these groups are conjugate to $\bar Q=R(Q)$ where $Q$ is the quaternion group of order $8$ generated by $\ii,\jj$. The intersection $T\cap Q$ is generated by $\ii$, a cyclic group of order $4$. Therefore $\bar\aA$ consists of the trivial subgroup, $\bar T$, $\bar Q$, and $\bar T\cap \bar Q \cong \ZZ/2$.
We describe the normalizers in $\SU(2)$. We have seen that $N(T)=\Span{T,\jj}$. One can verify that 
\[N(Q)=\left\langle \frac{\one+\ii}{\sqrt{2}},\frac{\one+\jj}{\sqrt{2}}\right\rangle\;\;\text{ and }\;\; N(T\cap Q) =N(T).\]
$N(Q)$ is the signed symmetric group on $3$ letters. The normalizers of $\bar T$, $\bar Q$, and $\bar T\cap \bar Q$ can be obtained from these by modding out by the centre $Z$. Then the object space of $\aA$ is given by
\[
\SO(3)/N(\bar T)_+ \bigvee \SO(3)/N(\bar Q)_+\bigvee \SO(3)/N(\bar Q \cap \bar T) _+
\]
and the corresponding quotient category  fits in a push-out diagram of categories
\begin{equation}\label{diag:copies}
\begin{tikzcd}
{[1]} \arrow[r,"d^2"] \arrow[d,"d^2"'] & {[2]}\arrow{d}\\
{[2]} \arrow{r} & \bar\aA
\end{tikzcd}
\end{equation}
We can think of the resulting category obtained by identifying two copies of $[2]$ along the edge $0\rightarrow 1$. 
Let $X:\aA \rightarrow\catTop$ be a left $\aA$-module. We will denote by $X^+$ and $X^-$ the restriction of the functor $\rho_* X$ to one of the copies of $d[2]^\op$. If $X=B(\ZZ,-)$ then one of them, say $X^+$, corresponds to $B(\ZZ,\SO(3))_1$. In fact the homotopy colimit diagram for this part simplifies to a homotopy push-out 
 $$
\begin{tikzcd}
\SO(3)/N(\bar T) \arrow{r}\arrow{d} & \SO(3)\times_{N(\bar T)}B\bar T \arrow{d}\\
\ast \arrow{r} & B(\ZZ,\SO(3))_1
\end{tikzcd}
$$
as in \cite{AGV17}.

Let us look at the space of $n$-simplices of $B(\ZZ,\SO(3))$. By Remark \ref{HomSpaceDec} we have a homotopy colimit decomposition for each of the functors $B_n(\ZZ,-)$, namely the functor of homomorphism spaces $\Hom(\ZZ^n,-)$. It is known that the connected components of $\Hom(\ZZ^n,\SO(3))$ other than the component $\Hom(\ZZ^n,\SO(3))_1$ of the trivial representation are all homeomorphic to $\SU(2)/Q$ \cite{TS08}. Then the cofibre $C_n$ of the inclusion map
$
\Hom(\ZZ^n,\SO(3))_1 \rightarrow \Hom(\ZZ^n,\SO(3)) 
$ 
fits in another cofibre sequence
$$
SU(2)/Q \rightarrow SU(2)/Q \times \pi_0 \Hom(\ZZ^n,\SO(3)) \rightarrow C_n
$$
where the inclusion is at the base point of the homomorphism space.
The number of connected components is computed in \cite{TS08} but we need a more explicit description.
We can compute the number of connected components from the homotopy decomposition of $\Hom(\ZZ^n,-)$.

\begin{lem}\label{pushoutSO}There is a push-out diagram
$$
\begin{tikzcd}
\Hom(\ZZ^n, \ZZ/2) \arrow{r} \arrow{d} & \Hom(\ZZ^n,(\ZZ/2)^2)/\Sigma_3 \arrow{d} \\
\ast \arrow{r} & \pi_0 \Hom(\ZZ^n, \SO(3))
\end{tikzcd}
$$
where $\Sigma_3$ permutes the non-trivial elements of $(\ZZ/2)^2$, and the top map is induced by one of the inclusions $\ZZ/2\subset (\ZZ/2)^2$.
\end{lem}
\begin{proof}
The decomposition corresponding to $X=\Hom(\ZZ^n,-)$ has two parts $X^+$ and $X^-$ where $X^+$ corresponds to the connected component $\Hom(\ZZ^n,-)_1$. Let $X^0$ denote the restriction  of $\rho_*X$ to the category $d[1]^\op$ obtained from the intersection of the two copies of $[2]$   in (\ref{diag:copies}). There is a push-out diagram of sets
$$
\begin{tikzcd}
\pi_0\,\hocolim{}X^0 \arrow{r} \arrow{d} & \pi_0\,\hocolim{}X^+ \arrow{d}  \\
\pi_0\,\hocolim{}X^- \arrow{r} & \pi_0\, \Hom(\ZZ^n, \SO(3)) 
\end{tikzcd}
$$
Since the $X^+$ part is connected the image of the $X^0$ part in the $X^-$ part will be identified to a single connected component. $X^-$ is a diagram of spaces on $d[2]^\op$ which sends a simplex $\sigma$ to the space $G\times_{N(\sigma)} X^-(\sigma_0)$. Here $0\rightarrow 1 \rightarrow 2$ is identified with the inclusion of subgroups $1\rightarrow \bar T\cap \bar Q \rightarrow \bar Q$. But we can also use the smaller diagram $1\rightarrow \bar Q$ to compute the homotopy colimit of $X^-$. From the latter we observe that the set of connected components is given by $\Hom(\ZZ^n,\bar Q)/N(\bar Q)$.
Similarly $X^+$ and $X^0$ can be described as homotopy push-out diagrams.  
The $X^0$ part corresponds to the image of the inclusion 
$$
\Hom(\ZZ^n,\bar T\cap \bar Q) \rightarrow \Hom(\ZZ^n,\bar Q)/N(\bar Q)
$$ 
and it is identified to a single point.
\end{proof} 

So far we have dealt with a fixed simplicial degree. Geometric realization commutes with push-outs and cofibre sequences. The cofibre of $B(\ZZ,\SO(3))_1\rightarrow B(\ZZ,\SO(3))$ is the geometric realization of the cofibres $C_\dt$. The following will be useful to understand the geometric realization of the push-out diagram in Lemma \ref{pushoutSO}.

\Lem{\label{quotientSO} There is a homotopy equivalence $(B(\ZZ/2)^2/\Sigma_3)^\wedge_\ell\simeq \ast$ for $\ell>2$. }
\Proof{Let us take the quotient in two steps. First take the quotient by the normal subgroup $\ZZ/3\subset \Sigma_3$, and then by $\Sigma_3/(\ZZ/3)\cong \ZZ/2$.
The action of $\ZZ/3$ on the $n$-simplices of $B(\ZZ/2)^2$ is free except at the degenerate simplex of the base point. Therefore there is a cofibre sequence
$$
B\ZZ/3 \stackrel{\sim_\ell}{\longrightarrow} E\ZZ/3 \times_{\ZZ/3} B(\ZZ/2)^2 \rightarrow B(\ZZ/2)^2/(\ZZ/3)
$$
where the first map is a mod-$\ell$ equivalence for $\ell>2$.
Thus $B(\ZZ/2)^2/(\ZZ/3)$ is mod-$\ell$ equivalent to a point and taking a further quotient by $\ZZ/2$ does not change the mod-$\ell$ homotopy type.
}

\Thm{\label{ModlEquivn=1}The homomorphism $R:\SU(2)\rightarrow \SO(3)$ induces a mod-$\ell$ equivalence
$$
B(\ZZ,\SU(2)) \rightarrow B(\ZZ,\SO(3))
$$
for all $\ell>2$.}
\begin{proof} Let $Z$ denote the centre of $\SU(2)$. Its classifying space $BZ$ is a topological group which acts freely on $B(\ZZ,\SU(2))$ and $B(\ZZ,\SO(3))_1$ is the orbit space of this action. Therefore $R$ induces a fibration $B(\ZZ,\SU(2))\rightarrow B(\ZZ,\SO(3))_1$. The fibre is homotopy equivalent to $B\ZZ/2$ hence it induces a mod-$\ell$ equivalence. We will show that the inclusion $B(\ZZ,\SO(3))_1\rightarrow B(\ZZ,\SO(3))$ is a mod-$\ell$ equivalence. Recall that at the level of $n$-simplices the cofibre is $C_n$. Since geometric realization commutes with cofibre sequences we have a cofibre sequence
$$
\SU(2)/Q \rightarrow \SU(2)/Q \times |X_\dt| \rightarrow |C_\dt|
$$
where $X_n=\pi_0 \Hom(\ZZ^n, \SO(3))$. By Lemma \ref{pushoutSO} we have a push-out diagram
$$
\begin{tikzcd}
{B\ZZ/2} \arrow{r} \arrow{d} & {B(\ZZ/2)^2/\Sigma_3} \arrow{d} \\
\ast \arrow{r} & {|X_\dt|}
\end{tikzcd}
$$
 Lemma \ref{quotientSO} implies that $|X_\dt|^\wedge_\ell \simeq \ast$, and consequently $|C_\dt|^\wedge_\ell$ is contractible. 
\end{proof}

Using the cohomology calculation in \cite{AGV17} we obtain

\Cor{Let $l$ be an odd prime. Then there is an isomorphism of graded rings
$$
H^*(B(\ZZ,\SO(3)),\ZZ/\ell) = \FF_\ell [p,y]/(y^2) 
$$
where both generators are of degree $4$. 
}

\section{Mapping spaces and spaces of homomorphisms} 

Given a topological group $G$ and a discrete group $\pi$, one may construct a ``mapping groupoid'' $\uHom(\pi, G)$. One
would like to know to what extent this formation of mapping object commutes with taking classifying spaces. In an ideal
situation, one would have
\begin{equation}
 \label{eq:3}
 B \uHom(\pi, G) \simeq \Map(B \pi, BG),
\end{equation}
the equivalence being taken over the space $BG$.

This holds when $\pi$ is $\ZZ$, since it is a rephrasing of the observation of D.~Sullivan that $EG \times_G G^{ad}
\simeq \Map( S^1, BG)$. In the case where $\pi$ is a finite $\ell$-group and $G$ is compact Lie, one has a mod-$\ell$ equivalence. This is Theorem \ref{mapping} below, due to Dwyer and Zabrodsky \cite{DZ87}.
If we take fibres over the maps to $BG$ in the putative equivalence \eqref{eq:3}, we obtain a comparison between the space
$\Hom(\pi, G)$ and the space of based maps $\Map_\ast (B \pi, BG)$. When $\pi$ is $\ZZ$, one has the well-known
equivalence $G \to \Omega B G$. When $\pi$ and $\pi_0(G)$ are finite $\ell$-groups, one has a mod-$\ell$ equivalence.

\subsection{Topologies on spaces of homomorphisms}

For any two groups $\pi$ and $G$, we write $\Hom(\pi, G)$ for the set of homomorphisms $\rho:\pi \to G$, and we write
$\Rep(\pi, G)$ for the quotient set by the conjugation action of $G$. The class of $\rho$ in $\Rep(\pi, G)$ will be
written $\langle \rho \rangle$, so $\langle \rho \rangle = \langle \rho' \rangle$ if and only if there exists some $g
\in G$ such that $g\rho(x)g^{-1} = \rho'(x)$ for all $x \in \pi$. 

Recall that the homomorphism space $\Hom(\ZZ^n,G)$ is topologized as a subspace of $G^n$. 

Suppose $\pi$ is a finite group generated by a finite set $S$ of elements. Then we may topologize $\Hom(\pi, G)$ as a subspace
of $G^S = \prod_{s \in S} G$. It is the case that the induced topology on $\Hom(\pi, G)$ does not depend on the
particular choice of generating set. The action of $G$
on $\Hom(\pi, G)$ is continuous.

\Pro{\label{Rep_disc} Let $\pi$ be a finite group and $G$ a compact Lie group. Give $\Hom(\pi, G)$ the topology specified above. Then the
quotient space $\Rep(\pi, G)$ is a discrete finite set.}
\Proof{ We will show that the space $\Rep(\pi, G)$ is finite. This will suffice to prove the statement, since $\Hom(\pi,
 G)$ is a Hausdorff space and $G$ is a compact Hausdorff group, so that the quotient set will inherit a Hausdorff
 topology. Being finite, it must be discrete.

 There exists an injective map of spaces
 \[ \Hom(\pi,G) \to \prod_{\alpha\in \pi}\Hom(\ZZ/|\alpha|,G )\]
 which is compatible with the conjugation action by $G$ and which remains injective after passing to $G$-quotients. It
 therefore suffices to show that $\Hom(\ZZ/m,G)/G$ is finite. For this embed $G$ into a unitary group $U(n)$. Let $S$
 denote the subgroup of elements of order dividing $\ell^k$ in a fixed maximal torus of the unitary group. Then $S$ has
 finite order, and the set of subgroups of $G$ conjugate under $U(n)$ to a cyclic subgroup of $S$ falls into finitely
 many conjugacy classes in $G$ by Quillen's argument in \cite[Lemma 6.3]{Q71}. This concludes the argument. }

\subsection{Maps between classifying spaces}

Suppose $\pi$ is a finite group and $G$ is a Lie group. One can define a topological groupoid of morphisms
$\uHom(\pi, G)$ by setting the set of objects to be the set of homomorphisms $\rho : \pi \to G$ and by declaring the
morphisms $\rho \to \rho'$ to consist of the subset of elements $g \in G$ such that $g \rho(x) g^{-1} = \rho'(x)$ for
all $x \in \pi$. It is easily verified that the set $\hom(\rho, \rho')$ of morphisms is either empty or consists of a coset $gZ_G(\rho)$ of
the centralizer of the image of $\rho$ in $G$. This centralizer, $Z_G(\rho)$, is a closed subgroup of $G$. We endow
$\hom(\rho, \rho')$ with the subspace topology, and this endows $\uHom(\pi, G)$ with the structure of a topological
category. We remark that the topology on the set of objects is discrete in this case---often when the term ``topological
category'' is used in the literature, this is the kind of category that is meant.

The topological groupoid $\uHom(\pi, G)$ constructed above is a special case of the internal mapping object of topological
categories; here the categories in question are $\pi$ and $G$, both viewed as topological categories having one
object. Consequently, there is a counit or `evaluation' functor
\[ \uHom(\pi, G) \times \pi \to G. \]
This counit functor admits a direct description as follows. A morphism in $\uHom(\pi, G) \times \pi$ consists of the data
of a natural transformation $\nu: \phi \to \psi$ and an element $p$ of $\pi$. The transformation $\nu: \phi \to \psi$
between two homomorphisms $\pi \to G$ is given by conjugation by an element $g \in G$ so that $g \phi = \psi g$, and an
element $p$ of $\pi$. The image of this under the counit map is $g \phi(p)$, alternatively described as $\psi(p)g$.

The induced map on classifying spaces of topological categories is
$B\uHom(\pi, G) \times B\pi \to BG$ and this has an adjoint, namely a canonical map
\[ B \uHom(\pi, G) \to \Map(B\pi, BG).\]

\Thm{\label{mapping}\cite{DZ87}
Let $\ell$ be a prime, $\pi$ be a finite $\ell$--group, and $G$ be a compact Lie group. Then the natural map
\[B\,\underline{\Hom}(\pi,G) \rightarrow \Map(B\pi,BG)\]
is a strong mod-$\ell$ equivalence.
} 

The following argument appears as part of the proof of the above, in \cite[Section 4]{DZ87}, and we set it aside for later reference.

\Pro{\label{pr:DZarg}
Let $\pi$ be a discrete group and let $G$ be a topological group, and let $f: \pi \to G$ be a homomorphism. If
$\Map(B\pi, BG)_{Bf}$ denotes the path component of $Bf$, with basepoint $Bf$, then there is a natural equivalence $\Omega \Map(B \pi, BG)_{Bf} \simeq G^{h \pi}$
where $\pi$ acts on $G$ via $\pi\cdot g = f(\pi) g f(\pi)^{-1}$.}

\Proof{
 By use of adjunctions, the space $\Omega \Map(B \pi, BG)_{Bf}$ may be identified with the subspace of the space of
 maps $\Map(B \pi, L BG)$ consisting of those maps yielding to $Bf \in \Map(B\pi, BG)$ after post-composition with the
 evaluation-at-the-basepoint map $LBG \to BG$. That is, we are studying the space of lifts 
 \[ \xymatrix{ & LBG \ar[d] \\ B \pi \ar[r] \ar@{-->}[ur] & BG } \]

 The space $LBG$ is a model for the Borel construction $EG \times_G G^{\text{ad}}$, and the space of lifts under
 consideration is therefore homeomorphic to the space of sections of the fibre bundle $E\pi \times_\pi G \to B\pi$ where the
 action of $\pi$ on $G$ is the composite of the adjoint action of $G$ on itself and the homomorphism $f$. The space of
 such sections is a model for $G^{h \pi}$. 
}

\Rem{\label{rem:subq} {\rm
Theorem \ref{mapping} establishes both the following facts. First, the path components of $\Map(B \pi, BG)$ are in bijection
with the set of representations $\Rep(\pi, G)$, these being the components of $B \uHom(\pi, G)$. Second, given a map
$\rho: \pi \to G$, one may take the component of $B\rho$ in the mapping space: $\Map(B\pi, B G)_{B \rho}$, which we give
the evident basepoint. Then the induced map on mapping spaces $\Omega B \uHom(\pi, G)_{B \rho} \to \Omega \Map(B \pi,
BG)_{B \rho}$ is the map $Z_G(\rho)= G^\pi \to G^{h \pi}$ from the centralizer of $\rho$ to the homotopy centralizer of
$\rho$ as defined in Proposition \ref{pr:DZarg}.
}}

\Def{\rm Given a space $X$ with a continuous action of a topological group $G$, we define the \textit{action groupoid}
 $X\agroupoid G$ to be the topological groupoid having $X$ as its space of objects and where the space of morphisms from $x$ to
 $y$ is the set $\{ g \in G \: : \: gx =y \}$, endowed with the topology of a subspace of $G$.}

At first glance, it appears that $\uHom(\pi, G)$ is the action groupoid $\Hom(\pi, G)\agroupoid G$ where $G$ acts on
homomorphisms by conjugation, but there is the question of the topologies chosen. Specifically, $\uHom(\pi, G)$ was
constructed by endowing the object set $\Hom(\pi,G)$ with the discrete topology, but unless $G$ is itself discrete, it
is unlikely that the action of $G$ on $\Hom(\pi, G)$ is continuous when $\Hom(\pi, G)$ is given the discrete topology.

With the choice of topology $\tau$, we can construct the action groupoid $\Hom(\pi, G) \agroupoid G$.

There exists a natural map of topological groupoids
\begin{equation}
 \label{eq:2}
 n: \uHom(\pi, G) \to \Hom(\pi, G) \agroupoid G
\end{equation}
given by the identity functions on the levels of sets.

\Pro{\label{Groupoid_Borel} Let $\pi$ be a finite group and $G$ a compact Lie group. The natural map
\[ Bn: B \uHom(\pi, G) \to B \Hom(\pi, G) \agroupoid G. \]
obtained by applying the classifying space functor to \eqref{eq:2} is a homotopy equivalence.}
\Proof{
 We first verify that the map induces a bijection $\pi_0(Bn)$ on sets of connected components. The map on object sets
 of the groupoids is the identity function $(\Hom(\pi, G), \text{disc}) \to (\Hom(\pi, G), \tau)$, and the induced map
 on $\pi_0$ is the induced map of equivalence classes $\langle \rho \rangle \mapsto \langle \rho \rangle$. The set
 $\pi_0(B \uHom(\pi, G))$ is $\Rep(\pi, G)$ with the discrete topology (by construction) whereas $\pi_0( B
 \Hom(\pi, G)\agroupoid G)$ is $\pi_0(\Rep(\pi, G))$ with the topology induced by $\tau$, but by Proposition \ref{Rep_disc} this is exactly the discrete
 topology as well, and therefore $\pi_0(Bn)$ is a bijection.

 It now suffices to prove that $Bn$ restricts to a homotopy equivalence in each connected component. To prove this, fix
 a homomorphism $\rho : \pi \to G$ and consider the components $B\uHom(\pi, G)_{\rho}$ and $(\Hom(\pi, G)\agroupoid G)_\rho$
 corresponding to $\rho$. Consider the full subgroupoids of $\uHom(\pi, G)_{\rho}$ and $(\Hom(\pi, G)\agroupoid G)_\rho$
 consisting of only the single object $\rho$. Each is the topological group $Z_G(\rho)$ and the map $n$ restricts to the identity
 map between these two groups. There is therefore a commutative diagram of topological groupoids
 \[ \xymatrix{ Z_G(\rho) \ar@{=}[r] \ar[d] & Z_G(\rho) \ar[d] \\ \uHom(\pi, G)_{\rho} \ar^n[r] & (\Hom(\pi, G)\agroupoid
 G)_\rho}. \]
 The two vertical maps are the inclusion of full subgroupoids on one object and therefore induce homotopy equivalences
 on classifying spaces. Applying $B$ to the diagram as a whole therefore yields the result.}

It is a folklore result that the classifying space of $\Hom(\pi, G) \agroupoid G$ is homeomorphic to the Borel
construction $EG \times_G \Hom(\pi, G)$.

 Let $X^\wedge_\ell$ denote the 
Bousfield-Kan $\ell$--completion of a space. A map $X\rightarrow Y$ is a mod-$\ell$ equivalence if and only if the
$\ell$--completed map $X^\wedge_\ell \rightarrow Y^\wedge_\ell$ is a homotopy equivalence, see Definition
\ref{def:lequiv}. For a compact Lie group $G$ and a finite $\ell$--group $\pi$ the map
\begin{equation}\label{map_space_completion}
\Map(B\pi,BG)^\wedge_\ell \rightarrow \Map(B\pi,BG^\wedge_\ell)
\end{equation}
induced by the completion map $BG\rightarrow BG^\wedge_\ell$ is a homotopy equivalence, \cite[Proposition 7.5]{BK02}. Note that the mapping space $\Map(B\pi,BG) $ is $\ell$--good since by Theorem \ref{mapping} it is mod-$\ell$ equivalent to a disjoint union of classifying spaces of centralizers. Each centralizer is a compact Lie group and the classifying space of a compact Lie group is $\ell$--good.

Let $\pi^\dt$ be a cosimplicial group. We say $\pi^\dt$ is a \textit{cosimplicial $\ell$--group} if $\pi^n$ is a finite $\ell$--group for all $n\geq 0$. We further require that $\pi^0$ is the trivial group $\set{1}$. 
This implies that the spaces of interest such as $\fat{\Hom(\pi^\dt,G)}$, $\fat{B\,\underline{\Hom}(\pi^\dt,G)}$ and $\fat{\Map(B\pi^\dt,BG)}$ are connected.

\Pro{\label{DZ_simplicial} 
Let $G$ be a compact Lie group, $\pi^\dt$ be a cosimplicial $\ell$--group. Then the natural map 
\[
\begin{tikzcd}
\fat{B\,\underline{\Hom}(\pi^\dt,G)} \rightarrow \fat{\Map(B\pi^\dt,BG^\wedge_\ell)} 
\end{tikzcd}
\]
is a mod-$\ell$ equivalence.
}
\Proof{
As a consequence of Theorem \ref{mapping} the natural map of simplicial spaces
\[
\phi^\dt: B\,\underline{\Hom}(\pi^\dt,G) \rightarrow \Map(B\pi^\dt,BG)
\]
is a degree-wise mod-$\ell$ equivalence. The homotopy equivalence in \ref{map_space_completion} implies that the completion map $BG\rightarrow BG^\wedge_\ell$ induces a mod-$\ell$ equivalence $\Map(B\pi^\dt,BG)\rightarrow \Map(B\pi^\dt,BG^\wedge_\ell)$ in each degree.
Then the fat realization of the composite map is also a mod-$\ell$ equivalence by a spectral sequence argument.}

\Thm{\label{hom_map}
Let $G$ be a compact Lie group, $\pi^\dt$ be a cosimplicial $\ell$--group.
There is a zigzag of mod-$\ell$ equivalences
\[ EG \times_G B(\pi^\dt, G) \leftarrow \fat{B\uHom(\pi^\dt, G) } \to \fat{\Map(B\pi^\dt,BG^\wedge_\ell)},\]
and this is natural in both the cosimplicial group $\pi^\dt$ and the compact Lie group $G$.
}
\Proof{ 
The simplicial space $X=\Hom(\pi^\dt,G)$ is good by Proposition \ref{good} and hence the quotient map $\fat{X}\rightarrow |X|$ is a homotopy equivalence.
Note that under the projection map an element $[x_n,u]$ in $\fat{X}$ maps to the equivalence class of the
 same representative in $|X|$ where $x_n$ is an $n$--simplex and $u$ is a point in the standard topological
 $n$--simplex. The $G$--action is on the first coordinate and hence the quotient is a $G$--map. When $X$ is good the
 homotopy inverse is defined in the same way, it is also a $G$--map. Since the homotopy inverse $|X|\rightarrow
 \fat{X}$ is a $G$--map and a homotopy equivalence the induced map $EG\times_G |X|\rightarrow EG\times_G \fat{X}$
between the Borel constructions is a homotopy equivalence. The fat realization and the Borel construction may both be
written as homotopy colimits, and therefore they commute with one another.
\[
EG\times_G |\Hom(\pi^\dt,G)| \simeq EG\times_G \fat{\Hom(\pi^\dt,G)} \simeq \fat{EG\times_G \Hom(\pi^\dt,G)}.
\]
Proposition \ref{Groupoid_Borel} gives a homotopy equivalence
\[
\fat{EG\times_G \Hom(\pi^\dt,G)} \overset{\sim}{\longleftarrow} \fat{B\,\underline{\Hom}(\pi^\dt,G)}
\]
since the fat realization respects homotopy equivalences.
Finally, using Proposition \ref{DZ_simplicial} finishes the proof.
}
 
Let $\Map_\ast(-,-)$ denote the pointed mapping space. 

\Cor{\label{hom_map_fib}
 Let $G$ be a compact Lie group, $\pi^\dt$ be a cosimplicial $\ell$--group (with $\pi^0=e$). Assume that $\pi_0(G)$ is an $\ell$--group. Then there is a natural zigzag of mod-$\ell$ equivalences
\[ 
B(\pi,G) \leftarrow \fat{ \Hom(\pi^\dt, G)} \rightarrow \fat{\Map_\ast(B\pi^\dt,BG^\wedge_\ell)} 
\]
} 
\Proof{
If we exploit the naturality of the zigzag in Theorem \ref{hom_map}, using the map from the trivial cosimplicial group
$\{e\}^\dt \to \pi^\dt$, we arrive at a diagram 
\[
\xymatrix{
EG\times_G B(\pi^\dt,G) \ar[d]& \fat{B \uHom(\pi^\dt, G) } \ar_\sim[l] \ar[d] \ar^{\sim_\ell}[r] & \fat{\Map(B\pi^\dt,BG^\wedge_\ell)}\ar[d] \\
BG \ar@{=}[r] & BG \ar^{\sim_\ell}[r] & BG^\wedge_\ell.
}
\]
The vertical maps appearing are fibrations and the induced natural map on fibres is the
zigzag
\[ B(\pi^\dt,G) \leftarrow \fat{ \Hom(\pi^\dt, G)} \rightarrow \fat{\Map_\ast(B\pi^\dt,BG^\wedge_\ell)}. \] 
The mod-$\ell$ homology of the fibres is finite dimensional in each degree.   For the space of homomorphisms this follows from the quotient map
$$
\Hom(\pi^n,G) \to \Rep(\pi^n,G)
$$
and for the mapping spaces from the fibration
$$
G\to \Map_*(B\pi^n,BG) \to \Map(B\pi^n,BG).
$$
Therefore when $\pi_0(G)$ is an $\ell$--group, the fibre lemma \cite[Ch II \S 5]{BK72} implies that these maps are mod-$\ell$ equivalences.
The result follows.
}

\Rem{\label{rem:counter}\rm{
The assumption on the connected components of $G$ is required. Corollary \ref{hom_map_fib} may fail in general. For example, let $G$ be the symmetric group $\Sigma_3$ on three letters. Inclusion of a Sylow $2$--subgroup $\ZZ/2 \rightarrow \Sigma_3$ induces a mod-$2$ equivalence between the classifying spaces. 
If Corollary \ref{hom_map_fib} were true in general it would imply that $B(\ZZ/2,\Sigma_3)$ is mod-$2$ equivalent to $B(\ZZ/2,\ZZ/2)=B\ZZ/2$ but this is not the case since 
$$
B(\ZZ/2,\Sigma_3) \simeq \vee^3 B\ZZ/2.
$$
}}

\Rem{\rm{ \label{component}
Let $\theta:\pi^\dt \rightarrow G$ be a fixed homomorphism. Under the assumptions of Corollary \ref{hom_map_fib} there is a natural zigzag of mod-$\ell$ equivalences
\[
B(\pi,G)_\theta\leftarrow \fat{ \Hom(\pi^\dt, G)_\theta} \rightarrow \fat{\Map_\ast(B\pi^\dt,BG)_{B\theta}}
\]
between the connected components at $\theta$.
}}
 
\section{Applications} 

We demonstrate two applications of the mod-$\ell$ (zigzag) equivalence between $B(\pi,G)$ and $\fat{\Map_\ast(B\pi^\dt,BG^\wedge_\ell)}$ as
proved in Corollary \ref{hom_map_fib}.

\subsection{Homotopy invariance} It is a well known fact that if a homomorphism $G\rightarrow H$ between groups is a
homotopy equivalence then the induced map $BG\rightarrow BH$ is also a homotopy equivalence. We should like to establish
a similar property for $B(\pi,G)$.

\Thm{\label{homotopy_inv} Let $G$ and $H$ be compact Lie groups, $\pi^\dt$ be a cosimplicial $\ell$--group. Assume that $\pi_0(G)$ and $\pi_0(H)$ are $\ell$--groups and
 there is a homotopy equivalence $BG^\wedge_\ell \rightarrow BH^\wedge_\ell$.
Then the induced map 
\[
\phi:B(\pi,G)^\wedge_\ell \rightarrow B(\pi,H)^\wedge_\ell
\]
is a homotopy equivalence.
}
\begin{proof} The induced map is given by the composite 
\[
\phi:B(\pi,G)^\wedge_\ell \rightarrow \fat{\Map_\ast(B\pi^\dt,BG^\wedge_\ell)} \rightarrow \fat{\Map_\ast(B\pi^\dt,BH^\wedge_\ell)} \rightarrow B(\pi,H)^\wedge_\ell
\] 
where the first and the last maps are homotopy equivalences provided by Corollary \ref{hom_map_fib}. 
\end{proof} 

Although taking a detour through finite $\ell$--groups is essential to the proof, using Corollary \ref{lie_modulo} we can extend
the result to
$B(\ZZ,G)$ using the commutative diagram
\[
\begin{tikzcd}
 B(\ZZ_\ell,G) \arrow{r} \arrow{d} & B(\ZZ,G) \arrow{d}\\
B(\ZZ_\ell,H) \arrow{r} & B(\ZZ,H)
\end{tikzcd}
\]

\Cor{\label{Z_hinv}Let $G$ and $H$ be compact Lie groups. Assume that $\pi_0(G)$ and $\pi_0(H)$ are $\ell$--groups and there is a homotopy equivalence $BG^\wedge_\ell \rightarrow BH^\wedge_\ell$. 
Then the induced map 
\[
B(\ZZ,G)^\wedge_\ell \rightarrow B(\ZZ,H)^\wedge_\ell
\]
is a homotopy equivalence.
}

\Ex{\rm{ In Theorem \ref{ModlEquivn=1} we proved that the homomorphism $\SU(2) \to \SO(3)$ induces a mod-$\ell$ equivalence $B(\ZZ,\SU(2)) \to B(\ZZ,\SO(3))$ when $\ell>2$. This is a special case of a more general result which is a consequence of Corollary \ref{Z_hinv}. Let $\Sp(n)$ denote the compact symplectic group. There is a homotopy equivalence
\[
B\Sp(n) ^\wedge_\ell \to B\SO(2n+1)^\wedge_\ell
\]
for every odd prime $\ell$ (see \cite{Fri75}). Theorem \ref{homotopy_inv} implies that for any cosimplicial $\ell$--group $\pi^\dt$ there is a homotopy equivalence
\[
B(\pi,\Sp(n))^\wedge_\ell \simeq B(\pi,\SO(2n+1))^\wedge_\ell
\]
and Corollary \ref{Z_hinv} gives a homotopy equivalence
\[
B(\ZZ,\Sp(n))^\wedge_\ell \simeq B(\ZZ,\SO(2n+1)) ^\wedge_\ell.
\]
}}

If $G$ is a compact Lie group such that $\pi_0(G)$ is an $\ell$-group, then $G$ furnishes an $\ell$--compact group,
\cite{DW94}. An $\ell$--compact group is a triple $(X,BX,e)$ where $BX$ is a pointed, connected, $\ell$--complete space,
$H^*(X, \ZZ/\ell)$ is finite, and $e: X\rightarrow \Omega BX$ is a homotopy equivalence. 

We can extend the definition of $B(\pi,-)$ to $\ell$--compact groups by using mapping spaces. 

\Def{\rm{\label{def2}
 For a cosimplicial group $\tau^\dt$ and an $\ell$--compact group $(X,BX,e)$ we define
\[
\bB (\tau,X) = \fat{\Map_\ast(B\tau^\dt, BX)}.
\]
 }}

 Let $(X,BX,e)$ and $(Y,BY,e')$ be $\ell$--compact groups. Then a homotopy equivalence $BX\rightarrow BY$ induces a homotopy equivalence
\[
\bB (\ZZ,X)^\wedge_\ell \rightarrow \bB (\ZZ,Y)^\wedge_\ell.
\]

\subsection{The Milnor--Friedlander conjecture} Let $k$ be an algebraically closed field, and $\ell$ be a prime distinct
from $p$, the characteristic of $k$. A reductive algebraic group $G$ over $k$ defines a sheaf on the \'etale site
$(\Sm|_k)_{\et}$ of smooth schemes over $k$. Applying the classifying space functor sectionwise gives a simplicial sheaf
$BG:(\Sm|_k)_{\et}\rightarrow \catsSet$. The functor of global sections, $\Gamma$, from simplicial sheaves to
simplicial sets is right adjoint to the constant-sheaf functor, $C$. The counit of this adjunction
$\epsilon: C\Gamma BG\rightarrow BG $ induces a map in \'{e}tale cohomology groups
\[
\epsilon^*:H^*_{et}(BG,\ZZ/\ell) \rightarrow H^*_{et}(C\Gamma BG, \ZZ/\ell)\cong H^*(BG(k),\ZZ/\ell)
\]
where $G(k)$ is the group of $k$--points, endowed with the discrete topology.
The Milnor--Friedlander conjecture asserts that $\epsilon^*$ is an isomorphism. This conjecture holds for $k=\bar{\FF}_p$ as proved by Mislin and Friedlander.

\Thm{\cite{FM84}\label{MF} 
Let $G(\CC)$ be a reductive complex Lie group, and let $G(\bar{\FF}_p)$ denote the $\bar{\FF}_p$--points of the associated Chevalley group $G_\ZZ$. 
Then there is a map 
\[
\Phi:BG(\bar{\FF}_p) \rightarrow B G(\CC),
\]
that is a mod-$\ell$ equivalence for every prime $\ell\not=p$.
}
 
We prove a similar theorem for $B(\ZZ,G(\CC))$.

Let $\WW$ denote the ring of Witt vectors of $\bar{\FF}_p$. This is a discrete valuation ring of characteristic 0 having
residue field $\bar{\FF}_p$. We denote the residue map by $r:\WW\rightarrow \bar{\FF}_p$. We fix an embedding
$s:\WW\rightarrow \CC$. We work over $S=\Spec\, \WW$ and denote by $G_S$ the base extension of the Chevalley group
$G_\ZZ$. Theorem \ref{MF} applies to generalized reductive group schemes $G_S$ in the sense of \cite[Section 2]{FM88}:
For such a group there is a map $\Phi:BG_S(\bar{\FF}_p) \rightarrow B G_S(\CC)$ that is a mod-$\ell$ equivalence for all
primes $\ell \neq p$. This map fits
in a homotopy commutative diagram
\[
\begin{tikzcd}
BG_S(\bar{\FF}_p) \arrow[dr,"\Phi_G"] & BG_S(\WW) \arrow{l}\arrow{d} \\
 & BG_S(\CC)
\end{tikzcd}
\]
see \cite[Remark 2.5]{FM86}. Here $G_S(\bar{\FF}_p)$ and $G_S(\WW) $ are discrete, and $G_S( \CC)$ is a Lie group.

\Pro{\label{MF_DZ}
Let $G_S$ be a reductive group over $S=\Spec\,\WW$. For a finite $\ell$--group $\pi$ the map $\Phi:BG_S(\bar{\FF}_p) \rightarrow B G_S(\CC)$ induces a mod-$\ell$ equivalence
\[
\tilde{\Phi}: \Map(B\pi, BG_S(\bar{\FF}_p)) \rightarrow \Map(B\pi, BG_S(\CC)) .
\]
}

We remark that the $\ell$-group $\pi$ is necessarily solvable.

\begin{proof}
Let us write $G=G_S$ for simplicity of notation.

We consider the zigzag
\begin{equation}
 \label{eq:4}
 BG (\bar \FF_p) \leftarrow BG(\WW) \to BG(\CC).
\end{equation}
Take mapping spaces $\Map(B\pi, \cdot)$ to obtain a zigzag
\begin{equation}
 \label{eq:5}
 \Map(B\pi, BG (\bar \FF_p) )\leftarrow \Map(B\pi, BG(\WW)) \to \Map(B\pi, BG(\CC)).
\end{equation}

We first claim that \eqref{eq:5} induces an isomorphism on connected components. In each case, the set of connected
components is the set of representations of $\pi$ 
\begin{equation}\label{pi0}
\Rep(\pi, G(\bar{\FF}_p) )\leftarrow \Rep(\pi, G(\WW) ) \rightarrow \Rep(\pi, G( \CC) )
\end{equation}
These maps are bijections by virtue of the main theorem of \cite{FM88}.

Consequently, it suffices to work component by component. Fix a homomorphism $\tilde \rho : \pi \to G(\WW)$ and
construct the diagram of homomorphisms
\begin{equation*}\label{diag_rho}
\begin{tikzcd}
 & \pi \arrow[ld,"\rho"'] \arrow[d,"\tilde\rho"] \arrow[rd,"\rho'"] & \\
G(\bar{\FF}_p) & G(\WW) \arrow{l}\arrow{r} & G(\CC)
\end{tikzcd}
\end{equation*}
and the maps between the associated components
\begin{equation}\label{diag_mapping}
\Map(B\pi, BG(\bar{\FF}_p) )_{B\rho} \leftarrow \Map(B\pi, BG(\WW) )_{B\tilde\rho} \rightarrow \Map(B\pi, BG( \CC)
)_{B\rho'}.
\end{equation}

We take the loop space of Diagram \eqref{diag_mapping}. 

There is an equivalence, natural in $H$, of the form
\[\Omega\Map(B\pi, BH)_{Bf} \to H^{h\pi}\]
(see Proposition \ref{pr:DZarg}) where $\pi$ acts by conjugation on $H$ via the homomorphism $f$. If $H$ is a discrete
group, such as $G(\bar{\FF}_p)$ or $G(\WW)$, the homotopy fixed point space is the fixed point space $H^\pi$,
\textit{viz} the centralizer $Z_{H}(f)$. If $H=K$ is a compact Lie group, then the Sullivan conjecture, in the form of
\cite[Theorem B(c)]{CarlssonEquivariantstablehomotopy1991a}, implies that the space of homotopy fixed points is
mod-$\ell$ equivalent to the fixed points, $(K^\wedge_\ell)^{h\pi} \simeq (K^\pi)^\wedge_\ell$.

We obtain a diagram
\begin{equation}
 \label{eq:7}
 \xymatrix{
Z_{G(\bar{\FF}_p)}(\rho) \ar^\sim[d] & Z_{G(\WW)}(\tilde{\rho}) \ar[l]\ar[r] \ar^\sim[d] & Z_{G(\CC)}(\rho') \ar[d] \\
G(\bar{\FF}_p)^{h\pi} & G(\WW)^{h\pi} \ar[l] \ar[r] & G(\CC)^{h\pi} 
} 
\end{equation}
where the action of $\pi$ is via the homomorphisms $\rho$, $\tilde\rho$, and $\rho'$ as appropriate.

Let $Z=Z_{G_S}(\tilde\rho)$ denote the centralizer of $\tilde{\rho}$ in $G_S$. The result in \cite[Theorem 4.2]{FM88} implies that
$Z$ is a  generalized reductive group over $S$. We have
\[
Z(\bar{\FF}_p)= Z_{G(\bar{\FF}_p)}(\rho)\quad \text{and}\quad Z(\CC)=Z_{G(\CC)}(\rho')
\]
since the fields are algebraically closed, \cite[pg. 28]{Jan87}. Let $K_Z$ denote a maximal compact subgroup of
$Z(\CC)$, and consider the subgroup generated by $K_Z$ and $\im(\rho')$ in $G(\CC)$. There is a surjective continuous
homomorphism from the compact group $K_Z \times \im(\rho')$ to this group, so it is itself compact. Let $K$ denote a
maximal compact subgroup of $G(\CC)$ containing both $\im(\rho')$ and $K_Z$. We remark that $K^\pi \subset Z(\CC) =
G(\CC)^\pi$ is a compact subgroup containing a maximal compact subgroup $K_Z$, so it follows that $K^\pi =
K_Z$. 

We now consider the commutative diagram 
\[
\begin{tikzcd} 
\Map(B\pi, BG(\bar{\FF}_p))_{B\rho} \arrow[r,"\tilde{\Phi}"] & \Map(B\pi,BG(\CC))_{B\rho'} & \Map(B\pi, BK)_{B
 \rho'} \arrow[l, "\sim"] \\
BZ(\bar{\FF}_p) \arrow[u,"\approx"] \arrow[r,"\Phi_Z"] & BZ(\CC) \arrow{u} & BK^\pi \arrow[l, "\sim"] \arrow[u, "\sim_\ell"]
\end{tikzcd}
\]
The right vertical arrow is a mod-$\ell$ equivalence by the homotopy equivalence in Theorem \ref{mapping} (see Remark
\ref{rem:subq}). We finish the proof by remarking that the Friedlander--Mislin theorem \cite{FM86} implies that $\Phi_Z$
is a mod-$\ell$ equivalence.
\end{proof}

\Thm{\label{main_repeat}Let $G(\CC)$ be a reductive complex Lie group and let $G(\bar{\FF}_p)$ denote the $\bar{\FF}_p$--points of the associated Chevalley group.
 Fix a prime $\ell$ other than the characteristic $p$ of $\FF_p$. 
 Then there is a zigzag 
\[
\varphi: EG(\bar{\FF}_p) \times_{G(\bar{\FF}_p)} B(\ZZ,G(\bar{\FF}_p)) \leftrightarrow EG(\CC) \times_{G(\CC)} B(\ZZ,G(\CC))
\]
of maps that are mod-$\ell$ equivalences.
}
The zigzag of maps depends on the embedding of the Witt vectors $\mathbb W \to \mathbb C$, but is natural in the group $G$.
\begin{proof}
Let $\pi^\dt$ denote
 the cosimplicial group $(\ZZ/\ell^k)^\dt$. Applying Proposition \ref{MF_DZ} in each degree gives a mod-$\ell$ equivalence
\[
\fat{\Map(B\pi^\dt,G(\bar{\FF}_p))} \to \fat{\Map(B\pi^\dt,G(\CC))} . 
\]
We can apply Theorem \ref{hom_map} to a maximal compact subgroup $K\subset G(\CC)$ to replace the target by $ EG(\CC) \times_{G(\CC)} B(\pi,G(\CC))$ up to mod-$\ell$ (zigzag) equivalence. The source, on the other hand, is homotopy equivalent to $EG(\bar{\FF}_p) \times_{G(\bar{\FF}_p)} B(\pi,G(\bar{\FF}_p)) $ since $G(\bar{\FF}_p) $ is a discrete group. Taking a (homotopy) colimit over $k$ and using Corollary \ref{lie_modulo} and Proposition \ref{loc_finite} finishes the proof.
\end{proof}

\begin{ex} \label{ex:noMF}
\rm{ We provide an example that illustrates Theorem \ref{main_repeat}, and also shows that the statement of the theorem fails without taking the Borel construction.

  We consider the reductive algebraic group $\GL_2$. Let $p$ be a prime. We start by describing certain abelian
  subgroups of the finite group $\GL_2(\FF_q)$ where $q$ is a power of $p$. Let $T(\FF_q)$ denote the subgroup
  consisting of diagonal matrices. The centre $Z(\FF_q)$ is isomorphic to the group of units $\FF_q^\times$. Matrices
  normalizing $T(\FF_q)$ either leave each element fixed or permute the diagonal. The normalizer group is denoted by
  $N(T(\FF_q))$ and is isomorphic to the semi-direct product $T(\FF_q)\rtimes \Sigma_2$ where $\Sigma_2$ acts by
  permuting the diagonal. In fact, $N(T(\FF_q))$ is the normalizer of any subgroup of $T(\FF_q)$ strictly containing the
  centre.

  Let $\ell$ be an odd prime distinct from $p$. There exists an integer $r\geq 1$ such that $q^r-1\equiv 0\mod\ell$. We
  may embed $\GL_2(\FF_q) \rightarrow \GL_2(\FF_{q^r})$, and since we are really interested in the colimit
  $\GL_2(\bar\FF_p)$, we may safely replace $\FF_q$ by the extension $\FF_{q^r}$. We therefore shall assume $q \equiv 1
  \pmod \ell$.

  Let $T_\ell(\FF_q)$ denote the $\ell$-torsion subgroup of the group $T(\FF_q)$ of order $(q-1)^2$. Since $\GL_2(\FF_q)$ has order $q(q-1)^2(q+1)$ this
  is a Sylow $\ell$-subgroup.

Let $\aA_\ell(\FF_q)$ denote the (discrete) poset of abelian $\ell$-subgroups of $\GL_2(\FF_q)$ that contain the centre.
The maximal elements of this poset are the conjugates of $T_\ell(\FF_q)$ and any two distinct conjugates intersect in
the $\ell$-torsion part of the centre---denoted $Z_\ell(\FF_q)$.

For the finite group $\GL(\FF_q)$ (or the locally finite group $\GL(\bar\FF_p)$) we can use  the discrete homotopy colimit decomposition as described in  \cite[\S 3]{O14}. By Proposition \ref{loc_finite}
$B(\ZZ,\GL_2(\FF_q))$ is mod-$\ell$ equivalent to $B(\ZZ_\ell,\GL_2(\FF_q))$ 
which    is homotopy equivalent to the homotopy colimit of
\begin{equation}\label{disc_hocolim}
\xymatrix{  & B( T_\ell(\FF_q)) \\
B( Z_\ell(\FF_q)) \ar[ur] \ar[dr] & \vdots\\
& B(T_\ell(\FF_q)) } 
\end{equation} 
If the $\ell$-adic valuation of $q-1$ is $s$, then $Z_\ell(\FF_q)$ is cyclic of order $\ell^s$ and $T_\ell(\FF_q)$ is a
product of two cyclic groups of order $\ell^s$. Moreover, it is possible to count the number of distinct conjugates of
$T_\ell(\FF_q)$ there are: namely $n_q = \frac{(q-1)^2 q (q+1)}{2 \ell^{2s}}$, the cardinality of
$\GL_2(\FF_q)/N(T_\ell(\FF_q))$.

For purposes of generalization, it is better to recast this homotopy colimit as the homotopy colimit of
\[ 
\begin{tikzcd}
\GL_2(\FF_q)/N(T_\ell(\FF_q)) \times BZ_\ell(\FF_q) \arrow{d}\arrow{r} &\GL_2(\FF_q)/T_\ell(\FF_q) \times_{\Sigma_2} BT_\ell(\FF_q) \arrow{d}\\
BZ_\ell(\FF_q) \arrow{r} & B(\ZZ,\GL_2(\FF_q))^\wedge_\ell .
\end{tikzcd}
\]
Here we write $\GL_2(\FF_q)/T_\ell(\FF_q) \times_{\Sigma_2} BT_\ell(\FF_q)$ in place of 
$\GL_2(\FF_q) \times_{N(T_\ell(\FF_q))} BT_\ell(\FF_q)$. The two spaces are homeomorphic, so this is not going to change
the homeomorphism type of the colimit, and $\GL_2(\FF_q)/T_\ell(\FF_q) \times_{\Sigma_2} BT_\ell(\FF_q)$ is preferred
because it indicates how $\GL_2(\FF_q)$ acts on the space. 
 
\medskip

We interested in taking a colimit over $q$ and determining $B(\ZZ,\GL_2(\bar \FF_p))$, at least up to mod-$\ell$
equivalence. To do this, let $\aA_\ell(\bar\FF_p)$ denote the poset of those abelian $\ell$-subgroups of $\GL_2(\bar\FF_p)$
containing the $\ell$-torsion part of the centre. This poset contains the subgroup $Z_\ell=\bigcup_q Z_\ell(\FF_q)$ and the
conjugates of $T_\ell = \bigcup_q T_\ell(\FF_q)$. Thus we have a homotopy push-out diagram
\[
\begin{tikzcd}
\GL_2(\bar\FF_p)/N(T_\ell) \times BZ_\ell \arrow{d}\arrow{r} &\GL_2(\bar\FF_p)/T_\ell \times_{\Sigma_2} BT_\ell \arrow{d}\\
BZ_\ell \arrow{r} & B(\ZZ,\GL_2(\bar\FF_p))^\wedge_\ell 
\end{tikzcd}
\]

The Borel construction of the conjugation action of $\GL_2(\bar\FF_p)$ on $B(\ZZ,\GL_2(\bar\FF_p))$ 
can be described from this diagram since the Borel construction commutes with homotopy colimits. 
One arrives at a homotopy pushout diagram
$$
\begin{tikzcd}
B N(T_\ell) \times BZ_\ell \arrow{d}\arrow{r} & EN(T_\ell)/T_\ell \times_{\Sigma_2} BT_\ell \arrow{d}\\
B\GL_2(\bar\FF_p) \times BZ_\ell \arrow{r} & X(\bar\FF_p) 
\end{tikzcd}
$$
for which the homotopy push-out is is mod-$\ell$ equivalent to $E\GL_2(\bar\FF_p) \times_{\GL_2(\bar\FF_p)} B(\ZZ,\GL_2(\bar\FF_p))$. 
\medskip

We compare this result to $\GL_2(\CC)$, or rather to the maximal compact subgroup $U(2)$. Let $\aA(\CC)$ denote the
topological category of intersections of maximal tori in $U(2)$. The objects in this poset consist of the centre $Z\cong U(1)$ and
the conjugates of a standard maximal torus $T\cong U(1)\times U(1)$.

From the homotopy decomposition \ref{disc_hocolim} we can show that the first mod-$\ell$ homology group of  $B(\ZZ,\GL_2(\bar\FF_p))$ is infinitely generated; whereas
the first homology of $B(\ZZ,U(2))_\ell$ vanishes   since $B(\ZZ,U(2))$ is simply connected, see \cite[Proposition 3.2]{AG15}. 

However, there is a similar homotopy push-out diagram 
 $$
\begin{tikzcd}
BN(T) \times BZ \arrow{r} \arrow{d} & EN(T )/T\times_{\Sigma_2} BT \arrow{d} \\
BU(2) \times BZ \arrow{r} & X(\CC) 
\end{tikzcd}
$$
where $X(\CC)$ is mod-$\ell$ equivalent to $EU(2)\times_{U(2)} B(\ZZ,U(2))$.
 $X(\bar\FF_p)$ is mod-$\ell$ equivalent to $X(\CC)$ 
since there is a natural mod-$\ell$ equivalence between the classifying spaces of $Z$, $T$, $N$, $U(2)$ and their discrete versions $Z_\ell$, $T_\ell$, $N_\ell$, $\GL_2(\bar\FF_p)$, respectively.

A similar argument applies to the algebraic group $\SL_2$.
}
\end{ex}


\end{document} 
